\documentclass[journal]{IEEEtran}
\usepackage{epigraph}
\usepackage{graphicx}
\usepackage{epstopdf}
\usepackage{numprint}
\usepackage{bm}
\usepackage{braket}
\usepackage[T1,T2A]{fontenc}
\usepackage[utf8]{inputenc}
\usepackage{hyperref}

\begin{document}
\title{Generalized Radon--Nikodym Spectral Approach.
  Application to Relaxation Dynamics Study.}
\author{Aleksandr Vasilievich Bobyl,
   Andrei Georgievich Zabrodskii,
   Mikhail Evgenievich Kompan,
   Vladislav Gennadievich Malyshkin,
   Olga Valentinovna Novikova,
   Ekaterina Evgenievna Terukova
   and Dmitry Valentinovich Agafonov.
   \thanks{
     A.B, A.Z, M.K and V.M are with
     Ioffe Institute, Politekhnicheskaya 26, St Petersburg, 194021, Russia, e-mails: bobyl@theory.ioffe.ru, Andrei.Zabrodskii@mail.ioffe.ru,
     kompan@mail.ioffe.ru, malyshki@ton.ioffe.ru ;
     O.N. is with Peter the Great St.Petersburg Polytechnic University,
     e-mail: novikova-olga1970@yandex.ru ; 
     E.T. is with R\&D Center TFTE LLC, St. Petersburg,
     e-mail:e.terukova@hevelsolar.com ;
     D.A. is with St. Petersburg State Technological Institute,
     e-mail: phti@lti-gti.ru ;    
 \$Id: RadonNikodymRelaxationDynamics.tex,v 1.67 2018/07/18 06:55:01 mal Exp $ $\$
   }
}

\date{November, 20, 2016}
\maketitle

\begin{abstract}
Radon--Nikodym approach to relaxation dynamics,
where probability density is built first
and then used to calculate observable dynamic characteristic is developed
and applied to relaxation type signals study.
In contrast with  $L^2$ norm approaches, such as
Fourier or least squares, this new approach
does not use a norm, the problem is
reduced to finding the spectrum of an operator (virtual Hamiltonian),
which is built in a way that eigenvalues represent
the dynamic characteristic of interest
and eigenvectors  represent probability density.
The problems of
interpolation (numerical estimation of  Radon--Nikodym derivatives is developed)
and  obtaining the distribution of relaxation rates 
from sampled timeserie are considered.
Application of the theory is demonstrated
on a number of model and experimentally measured
timeserie signals of degradation and relaxation processes.
Software product, implementing the theory
is developed.
\end{abstract}

\maketitle

\section{\label{intro}Introduction}
For the problem of relaxation signals study
the determination of dynamic characteristics
from experimentally collected data is critically important
in applications.
Initial data, typically collected in a form of timeseries,
is often analyzed applying such common techniques as
regressional analysis, Fourier/Laplace analysis, wavelet analysis,
Machine Learning, distribution extreme value analysis, among many others.
In this paper we developed, implemented in software
and demonstrated in applications
a new calculus-like technique\cite{malyshkin2015norm}, Radon--Nikodym derivatives (and their generalization),
to relaxation dynamics problems.
This technique is extremely effective
in relaxation signals study, where other
approaches (such as Fourier/Laplace analysis) are poorly applicable.
The synergy between Radon--Nikodym derivatives
and numerical analysis
(computer software  implementation is provided)
as well as approach effectiveness to relaxation processes study
are  the the most important features. Without apriori
assumptions (such us the number or different relaxation rates) 
the distribution of key dynamics characteristics of sampled signal
(such  as relaxation times, time derivative, etc.)
can be directly obtained from timeserie.
The approach
is reducible to matrix spectrum problem (without any $L^2$ norm involved),
and the distribution of dynamics characteristics
can be obtained from the distribution of eigenvalues in spectrum.
To the best of our knowledge this
is the only practical
approach allowing to obtain the distribution
of relaxation rates from sampled data without apriori assumptions.
Other approaches are not applicable for these reasons:
\begin{itemize}
\item Fourier/Wavelet transform  cannot be effectively applied
  to relaxation data study. If basis dimension
  is chosen high enough, then yes, all information
  is accumulated in components, but,
  in contrast with oscillatory type of signals,
  the distribution of relaxation rates cannot be obtained
   from the distribution of Fourier components.
  \item
    Laplace analysis is applicable only
     analytically, for sampled data discretization,
    rounding errors and sample finite size make Laplace approach inapplicable.
  \item ``Fitting the curve'' type of approaches
    can be effectively applied to the data such as in
    Figs. \ref{figLiIon10} and  \ref{figSCDischarge},
    but only with \textsl{apriori knowledge} about the number
    of different relaxation rates in sample (i.e. two in Fig. \ref{figLiIon10}
    and three in Fig. \ref{figSCDischarge}).
    Without this information, e.g. with multiple
    rates as in Fig. \ref{stockF},
    any ``fitting the curve'' approach fails.
  \item Re--sampling techniques
    (such as take initial $(x_l,f_l)$ timeserie and re-sample it to obtain derivative timeserie  $(x_l, (f(x_l)-f(x_l-\tau))/\tau)$
    typically fail for inability  to choose
    the parameter $\tau$.
    For real life data (such as presented in Fig. \ref{stockF}) a good  $\tau$ value does not exist.
  \item Additional important feature of our approach,
    is that it does not deploy any  $L^2$ norm
    (as Fourier or least squares do), it  uses operator spectrum instead.
    The approach is norm--free.
    For this reason it is not sensitive to sample outliers thus
    can be effectively applied to non--Gaussian samples,
    e.g. even with infinite second moment of dynamic characteristics
    in study.
\end{itemize}

Java written computer implementation is provided.
The code reads timeserie  data input $(x_l,f_l)$,
construct ``virtual quantum system Hamiltonian'',
corresponding to dynamic characteristics of interest
($f$, $df/dx$ or $(df/dx)/f$), diagonalize this ``Hamiltonian''
and output its spectrum. The  distribution
of dynamic characteristics of interest is described by
the distribution of obtained spectrum,
similar to the situation in random matrix theory\cite{guhr1998random}.
An application is demonstrated on a number of
practical examples, including
both model data and real life data
of degradation and relaxation dynamics problems.

\section{\label{RNinterpolatory}Radon--Nikodym Interpolation Problem}
Let us start with the simplest possible problem: interpolation problem.
Consider some timeserie ($x$ is time) of $M$ observations total
\begin{eqnarray}
  x_l&\to&f_l  \label{tsdata} \\
  l&=& 1\dots M \nonumber \\
  x_l&\le&x_{l+1} \nonumber
\end{eqnarray}
Interpolation problem is to estimate at given $x$ the
$f(x)$ and $df/dx$ using the data (\ref{tsdata}).
First we need some basis $Q_k(x)$, in this paper we will be using
polynomial bases. The results are invariant with respect
to basis choice (e.g. $Q_k(x)=x^k$ give identical result),
but numerical stability, especially in high dimension, depends strongly on basis choice\cite{malyshkin2015norm},
and the bases like Chebyshev or Legendre polynomials have to be
used in numerical calculations (see appendix \ref{appnum} for discussion).
Given the basis $Q_k(x)$ and observation weight function $\omega(x)$
(in this paper $\omega(x)=1$), 
obtain the moments:
\begin{eqnarray}
\Braket{Q_k}&=&\sum_{l=1}^{M} Q_k(x_l)  (x_l-x_{l-1})\omega(x_l) \label{qm} \\
\Braket{f Q_k}&=&\sum_{l=1}^{M} Q_k(x_l)  (x_l-x_{l-1})f_l \omega(x_l) \label{Cqm}\\
\Braket{df/dx Q_k}&=&\sum_{l=1}^{M} Q_k(x_l) (f_l-f_{l-1})\omega(x_l) \label{CDqm}
\end{eqnarray}
Given these moments Gramm matrix $G_{jk}=\Braket{Q_j Q_k}$
and the matrices $\Braket{f Q_j Q_k}$ and $\Braket{df/dx Q_j Q_k}$ can be readily obtained
either directly from (\ref{qm}), (\ref{Cqm}) and (\ref{CDqm}) sums
with $Q_j(x)Q_k(x)$ term, or, more computationally efficient,
using basis multiplication operator
(the $c_l^{jk}$ coefficients are available analytically for all practically interesting bases, see Appendix A of Ref. \cite{2015arXiv151005510G}
and references therein)
and first $2n-1$ moments $\Braket{Q_l}$, $\Braket{f Q_l}$ and $\Braket{df/dx Q_l}$;
$l=[0\dots 2n-2]$.
\begin{eqnarray}
  Q_j(x)Q_k(x)&=&\sum_{l=0}^{j+k} c_l^{jk}Q_l(x)
  \label{cmul}
\end{eqnarray}
Standard least squares interpolation $A_{LS}(y)$ is then:
\begin{eqnarray}
  A_{LS}(y)&=&\sum\limits_{i,j=0}^{n-1} Q_i(y)\left(G^{-1}\right)_{ij}\Braket{g Q_j} \label{leastsq}
\end{eqnarray}
here $G^{-1}$ is a matrix inverse to  Gramm matrix $\Braket{Q_j Q_k}$,
the $g$ is either $f$ of $df/dx$ depending on what need to be interpolated,
and $n$ is the dimension chosen.
Then the $A_{LS}(y)$ is $n-1$ order polynomial interpolating the function $g(y)$.
The (\ref{leastsq}) is direct interpolation of observable $g$,
obtained from  minimization of $L^2$ norm (cost function\cite{theodoridis2015machine})  
$\Braket{\left(g(x) -\sum_{i=0}^{n-1} \alpha_i Q_i(x)\right)^2}\to\min$
on $\alpha_i$,
what give $A_{LS}(y)=\sum_{i=0}^{n-1} \alpha_i Q_i(y)$
 in the form (\ref{leastsq})).

Radon--Nikodym interpolation $A_{RN}(y)$ is:
\begin{equation}
A_{RN}(y)=\frac{\sum\limits_{i,j,k,m=0}^{n-1}Q_i(y) \left(G^{-1}\right)_{ij} \Braket{g Q_jQ_k} \left(G^{-1}\right)_{km} Q_m(y)}
  {\sum\limits_{i,j=0}^{n-1}Q_i(y) \left(G^{-1}\right)_{ij}Q_j(y)}
  \label{RNsimple}
\end{equation}
The (\ref{RNsimple}) is a ratio of two polynomials
of  $2n-2$ order each, an estimator of stable form\cite{malha}.
In contrast with  $A_{LS}(y)$, which is obtained as one--stage
interpolation of observable $g$, the  $A_{RN}(y)$ is obtained in two stages:
\begin{itemize}
\item Obtain localized at $x=y$ probability density $\rho_y(x)=\psi^2_{y}(x)$,
  where
\begin{eqnarray}
  \psi_{y}(x)&=&\sum\limits_{i,j=0}^{n-1}  Q_i(x)  \left(G^{-1}\right)_{ij}  Q_j(y)
  \label{psiyX}
\end{eqnarray}
  is localized at $x=y$.
  The (\ref{psiyX}) is actually identical to least squares answer
  (\ref{leastsq}) with the replacement of  $\Braket{g Q_j}$ by $Q_j(y)$.
  Also note that $1/\psi_{x}(x)$ is equal to Christoffel function,
  which define Gauss quadratures weights\cite{totik}
  for the measure $d\mu=\omega(x)dx$.
\item  Average $g(x)$ with obtained probability density as
  \begin{equation}
    A_{RN}(y)=\frac{\Braket{\rho_y(x) g(x)}}{ \Braket{\rho_y(x)}}
    =
    \frac{\int   \psi_y^2(x) g(x)\omega(x) dx}{ \int   \psi_y^2(x)\omega(x) dx}
    \end{equation}
  to obtain  (\ref{RNsimple}) as a ratio of two quadratic forms,
  the ratio of two polynomials of  $2n-2$ order
  in case of polynomial basis.
\end{itemize}
The answer (\ref{RNsimple}) is plain Nevai operator\cite{nevai}
(with its property of absolute convergence to $g(x)$  with $n$ increase)
and can be viewed as numerical estimation of $\frac{d\mu}{d\nu}=\frac{g dx}{dx}$
considered as Radon--Nikodym derivative\cite{kolmogorovFA}.
 The (\ref{RNsimple}) answer is typically the most convenient
one among other available,
because it requires only one measure $d\nu$ to be positive. Other
answers\cite{BarrySimon,2015arXiv151005510G}
require both measures $d\mu$ and $d\nu$ to be positive.
This Radon--Nikodym interpolation (\ref{RNsimple}) has several
critically important advantages\cite{2015arXiv151005510G,2015arXiv151101887G,malyshkin2015norm}
compared to least squares
interpolation (\ref{leastsq}):
stability of interpolation, there is no divergence
outside of interpolation interval,
oscillations near interval edges are very much suppressed,
even in multi--dimensional case\cite{2015arXiv151101887G}.
These advantages come from the very fact,
that probability density is interpolated first,
then the result is obtained by averaging with this, always positive, interpolated
probability.

Another issue should be stressed here. For $df/dx$ interpolation
one should use the (\ref{CDqm}) moments for $g=df/dx$.
If one, instead of direct interpolation of $g=df/dx$,
 interpolate  $g=f$ (using (\ref{Cqm}) moments)
and then differentiate interpolating expression
(\ref{leastsq}) or (\ref{RNsimple})
the result will be incorrect\cite{2015arXiv151101887G}.

\subsection{\label{RNdem} Demonstration of Radon--Nikodym Interpolation in 1D and 2D cases}
Before we go further let us demonstrate the result of
Radon--Nikodym interpolation on standard data set\cite{2015arXiv151101887G} in 1D and 2D cases.

\begin{figure}[!t]
\includegraphics[width=4cm]{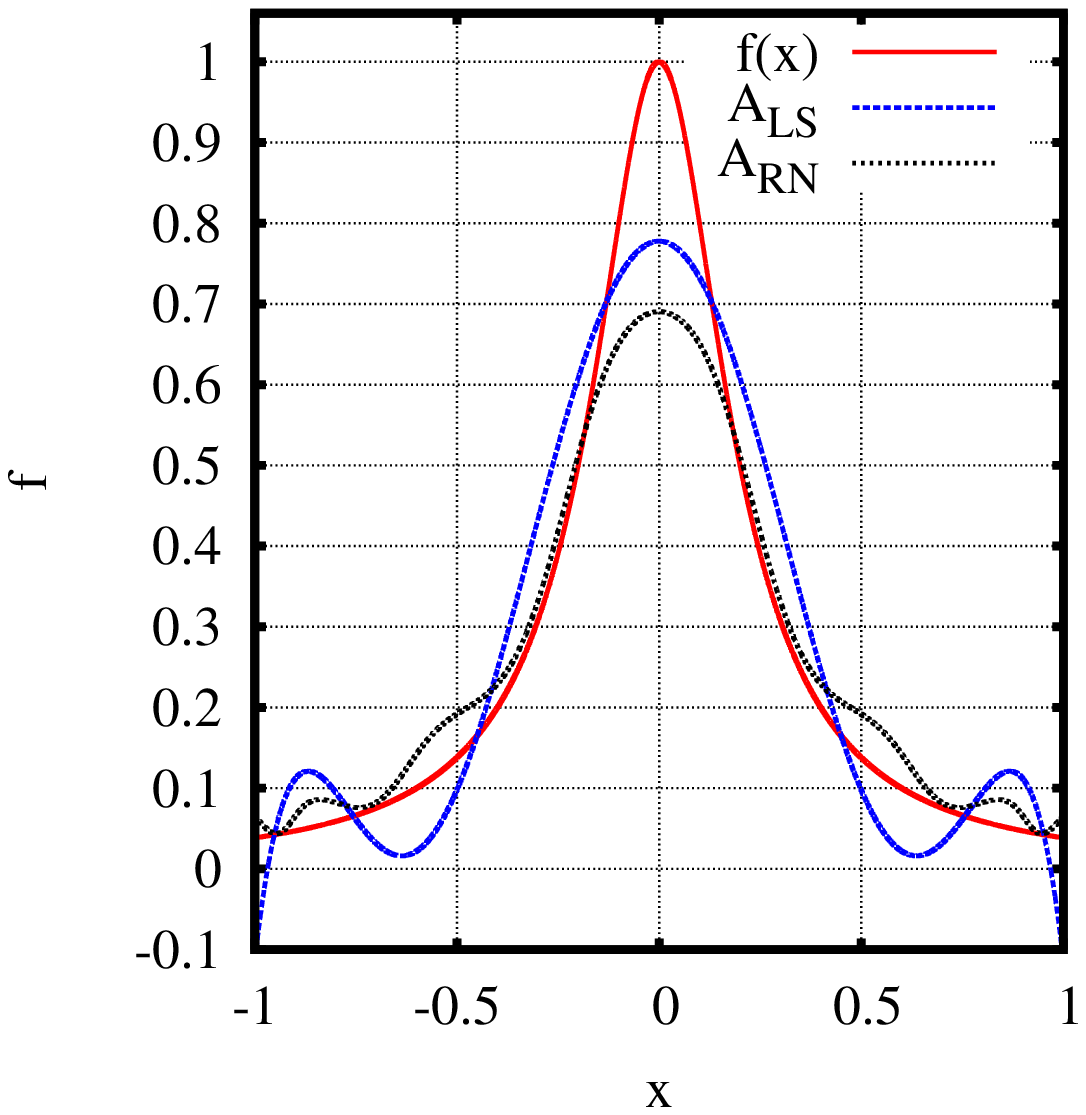}
\includegraphics[width=4cm]{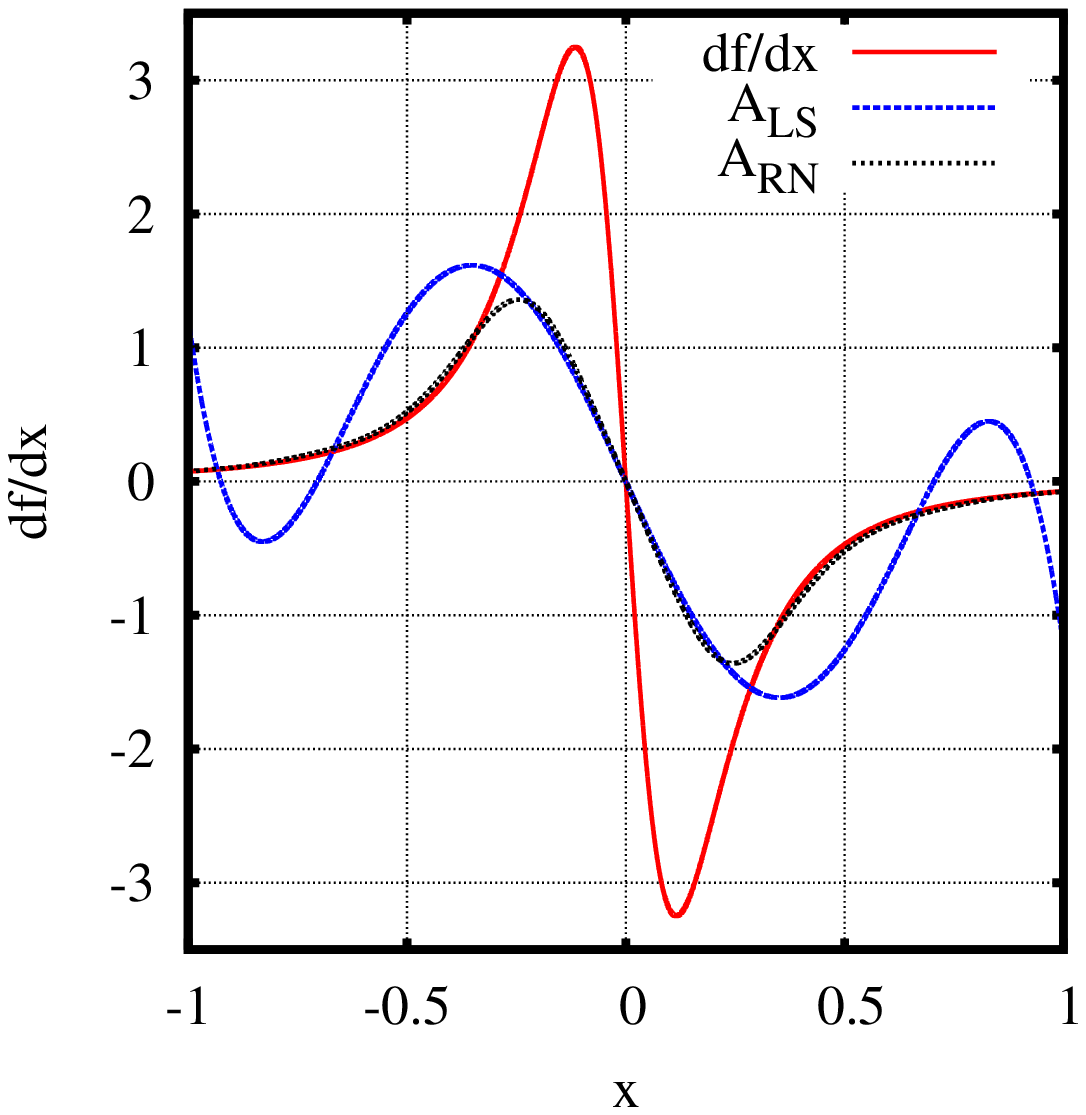} \\

\includegraphics[width=5.5cm]{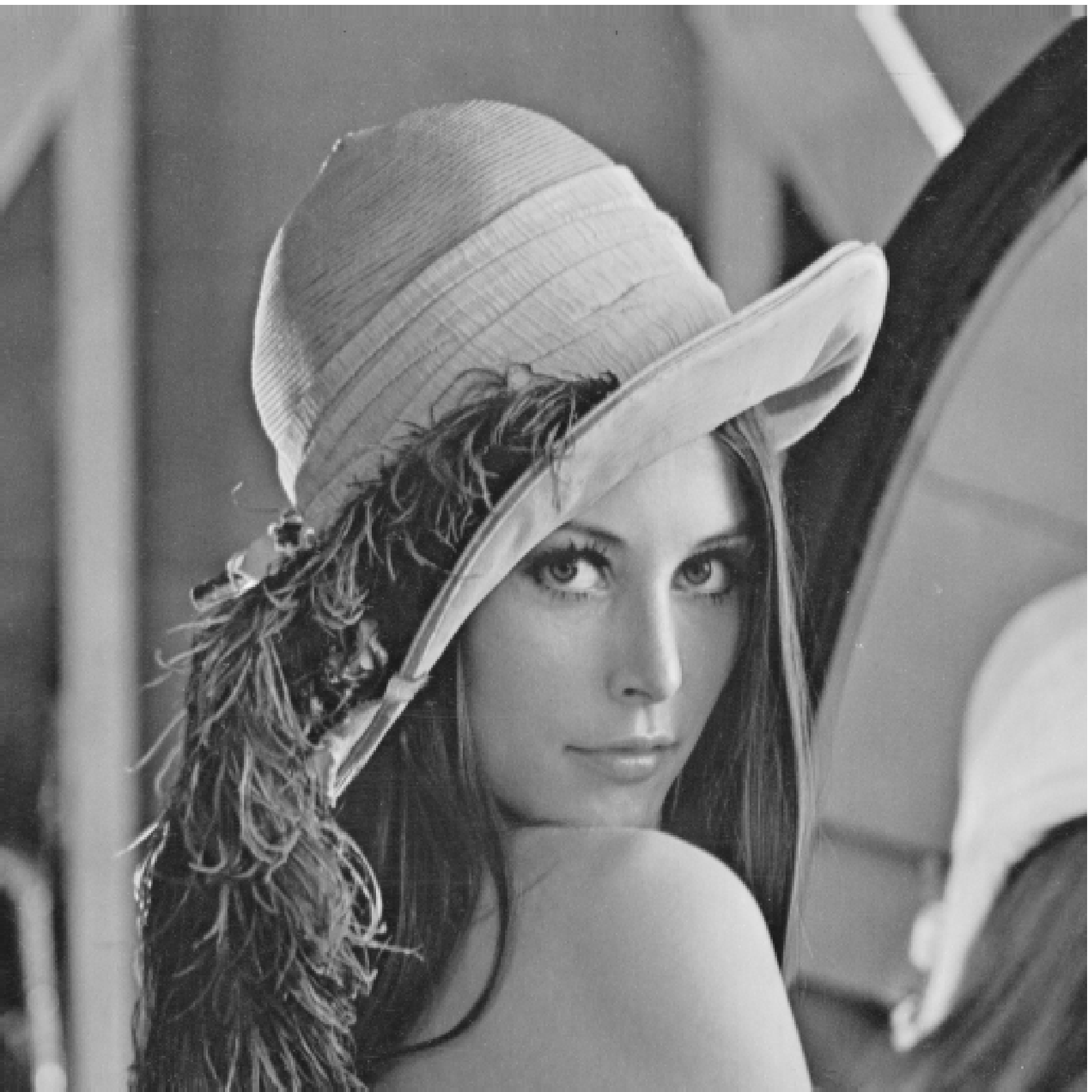}
\includegraphics[width=5.5cm]{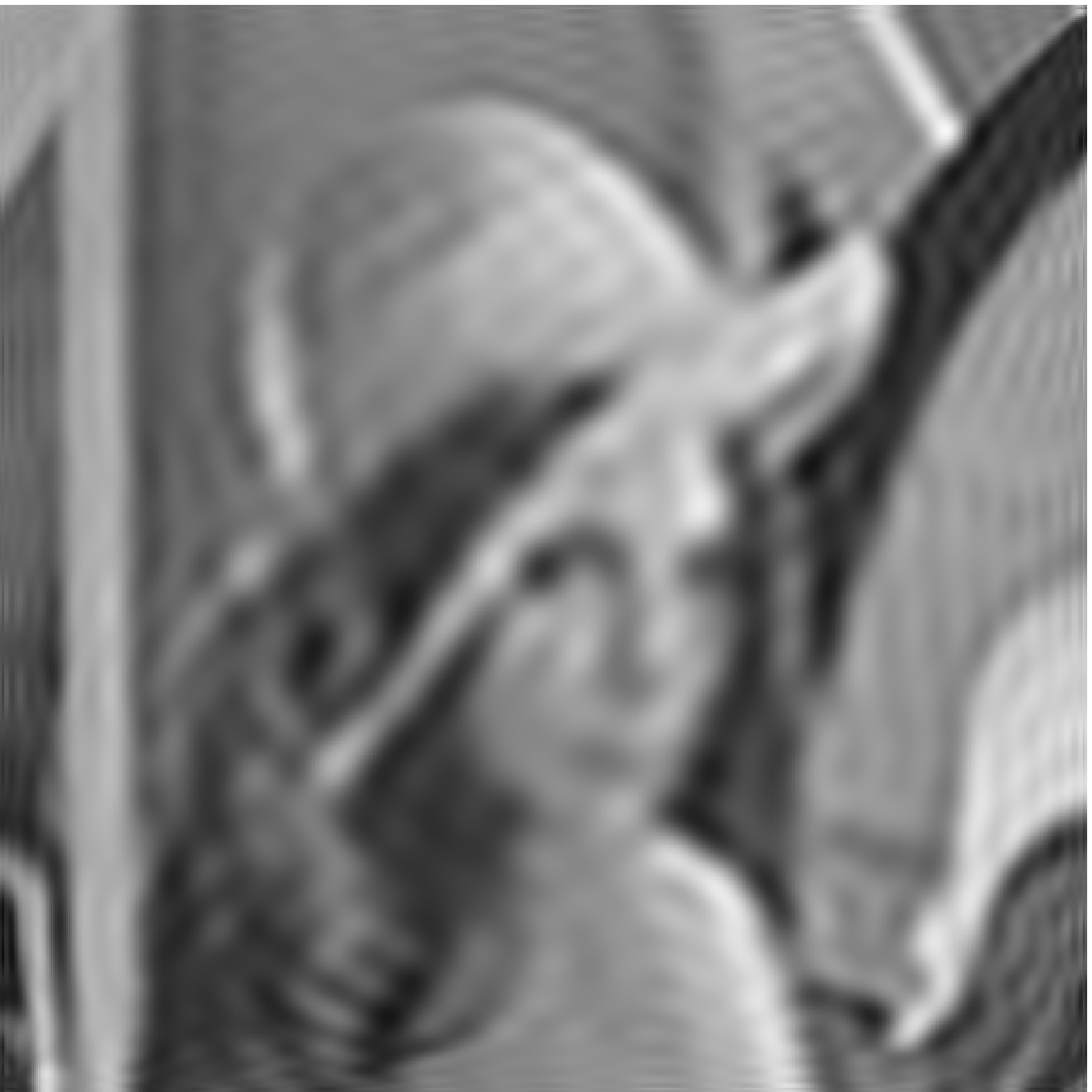}
\includegraphics[width=5.5cm]{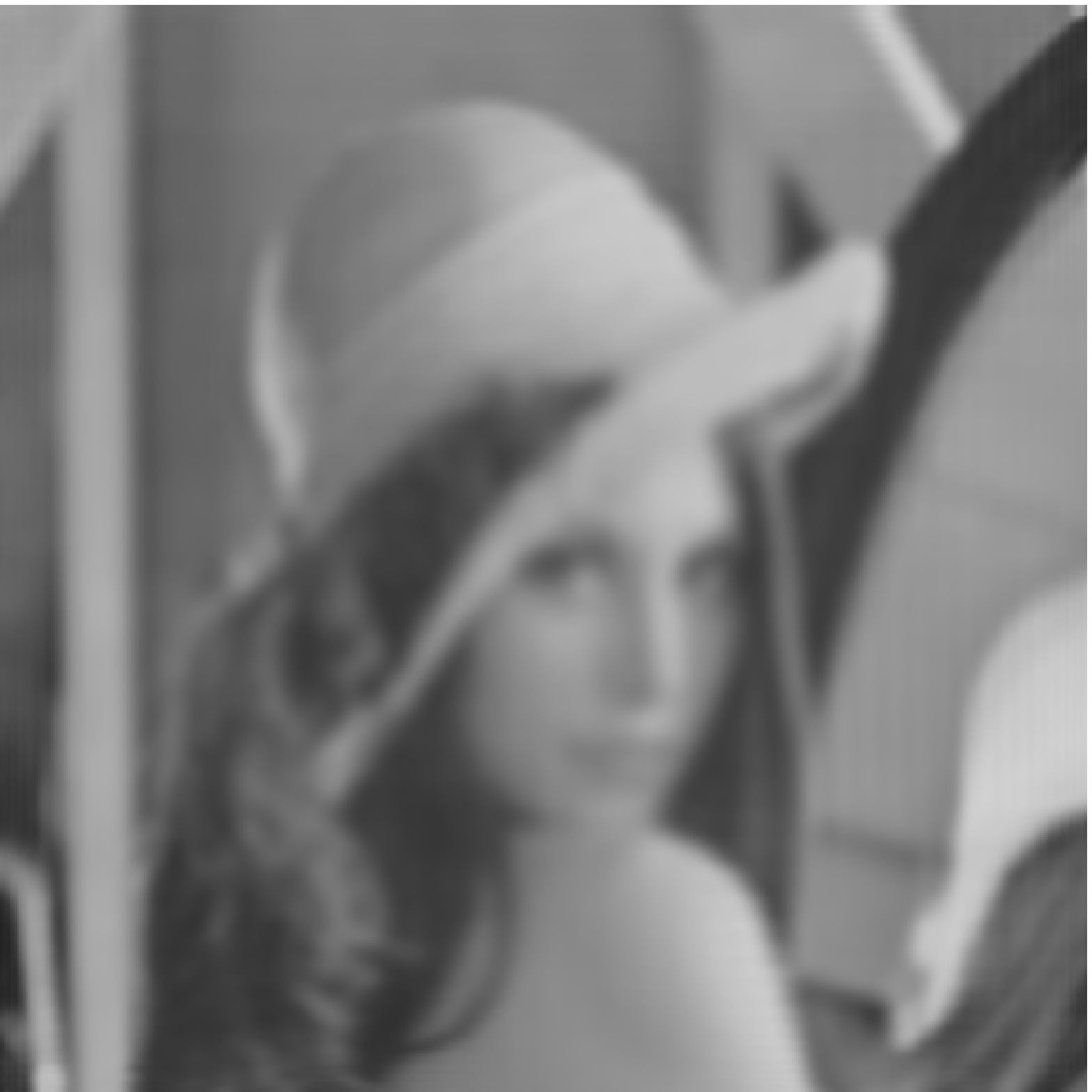}
\caption{\label{figinterpo}
  Top:Least squares $A_{LS}$ and Radon--Nikodym $A_{RN}$ approximations
  for Runge function(left) and  Runge function derivative (right)
  for $n=7$.
  Bottom:  Original Lenna image, least squares and Radon-Nikodym 
  for $n_x=n_y=50$.
}
\end{figure}

Consider (often used for algorithms testing) Runge function
\begin{equation}
  f(x)=\frac{1}{1+25x^2} \label{rungeF}
\end{equation}
with the measures
\begin{eqnarray}
  \Braket{f Q_{k}}&=& \int_{-1}^{1} f(x) Q_k(x) dx
  \label{mom1D} \\
  \Braket{df/dx Q_{k}}&=& \int_{-1}^{1} \frac{d f(x)}{dx} Q_k(x) dx
  \label{momd1D}
\end{eqnarray}
Then apply (\ref{leastsq}) and (\ref{RNsimple}) interpolating formulas.
The results are presented in Fig.\ref{figinterpo}.

In Fig. \ref{figinterpo} least squares and Radon--Nikodym
interpolations are presented for $n=7$ and the measures (\ref{mom1D}) and (\ref{momd1D}).
One can see that near edges oscillations are much
less severe, when  Radon--Nikodym approximation (polynomials ratio)
is used for interpolation.
This is because of the very fact that Radon--Nikodym approximation
has probability (not the observable value) interpolated,
and the result is obtained by averaging with this, always positive,
probability density.

Consider 2D case of image grayscale intensity interpolation
for Lenna 512x512 image\cite{wiki:lenna}
(presented in Fig. \ref{figinterpo}).
As a measure use the sum over image pixels $(t_x,t_y)$ and consider basis
index as double index $\bm{k}=(k_x,k_y)$:
\begin{eqnarray}
 && \Braket{f Q_{\bm{k}}}= \sum_{t_x=0}^{d_x-1}\sum_{t_y=0}^{d_y-1} f(t_x,t_y)
  Q_{k_x}(t_x)Q_{k_y}(t_y)  \label{momsample}
\end{eqnarray}
With a proper $Q_{\bm{k}}$ basis selection\cite{2015arXiv151101887G}
numerically stable results can be obtained
for up to $100\times 100$ elements in basis, i.e. for $10,000$ basis functions.
The least squares interpolation, same as in 1D case,
present typical for least squares intensity oscillations near image edges,
while  Radon--Nikodym has these oscillations very much suppressed.
Another important feature of Radon--Nikodym is that
it preserves the sign of interpolated function, i.e.
the grayscale intensity $f$ never become negative,
what may happen easily for least squares.

The difference between Radon--Nikodym and least squares
are of qualitative type
and can be summarized as following:
\begin{itemize}
\item For Radon--Nikodym (\ref{RNsimple})
   interpolated value at $y$ is
   sample data averaged with positive weight $d\mu_{\psi}=\psi_{y}^2(x)\omega(x)dx$ localized at $y$
   (\ref{psiyX}),
   i.e. it  is some weighted combination of sample observations.
   Least squares (\ref{leastsq})  interpolation
   does not correspond to a weighted combination of sample observations.
   This lead to interval edge oscillations
   being very much suppressed for Radon--Nikodym.
 \item Radon--Nikodym is a ratio of two quadratics forms
   built on basis functions. Least squares is a linear
   combination of basis functions.
   This lead to interpolation function sign preservation
   and no divergence at $y\to\infty$ for  Radon--Nikodym.
 \item Because Radon--Nikodym does not deploy any $L^2$ norm
   it is much less sensitive to data outliers.
   For least squares even a single data outlier
   can completely break the result.
 \item Radon--Nikodym requires more computational
   resources: $n\times n$ matrix $\Braket{g Q_jQ_k}$
   (obtained from first $2n-1$ input moments)
   vs. a vector of first $n$ input moments $\Braket{g Q_j}$ for least squares.
 \item Most important, Radon--Nikodym  can be generalized
   as an approach where learned knowledge
   is stored in matrix spectrum and then extracted using projection operators.
   This to be considered next.
\end{itemize}

\section{\label{RNQMbasis}Generalized Radon--Nikodym. Quantum Mechanical Spectral Description of Classical Experiment}
Considered in section \ref{RNinterpolatory} interpolation
problem is the first application of Radon--Nikodym approach.
The key element of the approach is that it based on matrices
$\Braket{Q_j Q_k}$ and $\Braket{g Q_j Q_k}$, not vectors $\Braket{g Q_j}$,
as in least squares (\ref{leastsq}).
The answers (\ref{leastsq}) and (\ref{RNsimple}) are written in
arbitrary polynomial basis $Q_i(x)$. However the matrix approach allow us to introduce
special, ``natural basis (\ref{psiC})'' $\psi^{[i]}(x)$, that diagonalize
the matrices $\Braket{Q_j Q_k}$ and $\Braket{g Q_j Q_k}$ simultaneously.
Consider Hamiltonian of a quantum system, the eigenfunctions of which
(\ref{psiC}),
are the solutions of generalized eigenvalues problem (\ref{GEV}), 
and diagonalize\cite{lapack} these two matrices simultaneously
$\Braket{\psi^{[i]}\psi^{[j]}}=\delta_{ij}$ and $\Braket{g \, \psi^{[i]}\psi^{[j]}}=\lambda^{[i]} \delta_{ij}$
\begin{eqnarray}
  && M^{L}_{jk}=\Braket{g \, Q_j Q_k} \\
  && M^{R}_{jk}=\Braket{ Q_j Q_k} \\  
&&\sum\limits_{k=0}^{n-1} M^{L}_{jk} \alpha^{[i]}_k =
  \lambda^{[i]} \sum\limits_{k=0}^{n-1} M^{R}_{jk} \alpha^{[i]}_k
  \label{GEV} \\ 
&&  \psi^{[i]}(x)=\sum\limits_{k=0}^{n-1} \alpha^{[i]}_k Q_k(x)
\label{psiC} \\
&&A_{RN}(y)=\frac{\sum\limits_{i=0}^{n-1} \lambda^{[i]} \left(\psi^{[i]}(y)\right)^2}
{\sum\limits_{i=0}^{n-1} \left(\psi^{[i]}(y)\right)^2}
\label{RNsimpleEVBasis} \\
&&A_{LS}(y)=\sum\limits_{i=0}^{n-1} \psi^{[i]}(y) \Braket{g \, \psi^{[i]}}
\nonumber
\\
&& =\sum\limits_{i=0}^{n-1} \lambda^{[i]} \psi^{[i]}(y)  \Braket{\psi^{[i]}}
 \label{leastsqEVBasis}
\end{eqnarray}
Probability states $\psi^{[i]}(x)$ from (\ref{psiC}) are normalized 
in a usual way: $\sum_{j,k=0}^{n-1} \alpha^{[s]}_j M^{R}_{jk} \alpha^{[t]}_k =\delta_{st}$
and
$\sum_{j,k=0}^{n-1} \alpha^{[s]}_j M^{L}_{jk} \alpha^{[t]}_k =\lambda^{[t]}\delta_{st}$.
The  (\ref{RNsimpleEVBasis}) and (\ref{leastsqEVBasis})
are {\bf identical} to  (\ref{RNsimple}) and  (\ref{leastsq}),
but written in the basis of  eigenfunctions $\psi^{[i]}(x)$
of  (\ref{GEV}) Hamiltonian.
This interpretation of Radon--Nikodym interpolation (\ref{RNsimple})
as obtaining two matrices from experimental observations,
introduction of (defined by two these matrices) Hamiltonian of a quantum system,
diagonalization (\ref{GEV}) and obtaining 
(\ref{RNsimpleEVBasis})
as an average with weights
equal to the square of projection (because
$\psi^{[i]}(y)=\Braket{\psi_{y}(x)\psi^{[i]}(x)}$ ;
also note that Christoffel function is $1/\sum\limits_{i=0}^{n-1} \left(\psi^{[i]}(y)\right)^2$)
 of a state  with specific $y$ (\ref{psiyX})
to eigenfunction (\ref{psiC})  is very fruitful.
The knowledge, obtained from the data,
is now accumulated in matrix spectrum,
and then can be extracted by projection operators.
This give a number of important advantages\cite{malyshkin2015norm}
with respect to standard approaches,
such as representation of learned knowledge in regression coefficients.

However, the most important generalization of the approach
is to go from localized $\psi_y(x)$ states (\ref{psiyX})
to arbitrary $\psi(x)$, e.g. the (\ref{psiC}) eigenvectors.
Then,
for example, minimal, maximal and average value of $g$ can be
easily estimated as 
$\min \lambda^{[i]}$, $\max \lambda^{[i]}$ and $\sum_{i=0}^{n-1}\lambda^{[i]}/n$.
This is drastically different from, say, regressional type of approach
where only function interpolation can be obtained from sample moments.
Obtained $\lambda^{[i]}$ spectrum
can be interpreted as the distribution of $g$ from sample,
what give a number of advantages in relaxation processes
study. For example\cite{ArxivMalyshkinMuse}, consider $n=2$,
then,  like a difference between median and average,
generalized
skewness  estimator of a $g(x)$ process
can be introduced:
\begin{eqnarray}
  \widetilde{\Gamma}&=&\frac{2\overline{g}-\lambda^{[\min]}-\lambda^{[\max]}}
            {\lambda^{[\min]}-\lambda^{[\max]}}
            \label{skewnesslikeG} \\
  \overline{g}&=&\frac{\Braket{gQ_0}}{\Braket{Q_0}}
  \label{gaver} 
\end{eqnarray}
Regular skewness estimator $\Braket{\left(g-\overline{g}\right)^3}$
requires 4 moments to calculate
$\Braket{1}$, $\Braket{g}$ 
$\Braket{g^2}$, $\Braket{g^3}$
and is not applicable to strongly non--Gaussian processes,
e.g. those with infinite $\Braket{g^2}$ or $\Braket{g^3}$.
The $\widetilde{\Gamma}$ 
skewness estimator requires
6 moments to calculate:
$\Braket{Q_0}$, $\Braket{Q_1}$, $\Braket{Q_2}$,
$\Braket{gQ_0}$, $\Braket{gQ_1}$, $\Braket{gQ_2}$,
(all of them are finite,
even when $\Braket{g^2}$ is infinite).
The $2\times 2$ matrices
$\Braket{gQ_jQ_k}$ and $\Braket{Q_jQ_k}$
can be readily
obtained from these moments, eigenvalues problem (\ref{GEV})
 solved by solving quadratic equation $0=\det\| \Braket{gQ_jQ_k}-\lambda\Braket{Q_jQ_k}\|$;
$\min g=\lambda^{[\min]}$, $\max g=\lambda^{[\max]}$
obtained, and $\widetilde{\Gamma}$ from (\ref{skewnesslikeG}) calculated.
When $Q_k(x)=Q_k(g)$
the $\widetilde{\Gamma}$
is proportional to regular skewness estimator. In this case
only 4 out of 6 moments required for $\widetilde{\Gamma}$ calculation
are independent. Another important result
of the approach is  a
``replacement to standard deviation'' as $\lambda^{[\max]}-\lambda^{[\min]}$,
that is finite even for the processes with infinite $\Braket{g^2}$.
This make the approach extremeny attractive to the
study of signals with spikes \cite{malyshkin2018spikes}.

\subsection{\label{RNoptimization}
  Quantum Mechanical Approach to
  Global Optimization Problem}
The Radon--Nikodym approach also give a  new look
to optimization problem\cite{2015arXiv151005510G}.
Instead of regualar  optimization problem
\begin{eqnarray}
  g(\bm{x})&\to&\min
  \label{gmin}
\end{eqnarray}
Consider optimization 
problem in ``quantum mechanics style''
(note argument and basis function index
can be considered multi-dimensional, e.g.
$Q_{{\bm{k}}}(\bm{x})=Q_{k_x}(x)Q_{k_y}(y)Q_{k_z}(z)\dots$):
\begin{eqnarray}
  &&\psi(\bm{x})=\sum_{\bm{k}}\alpha_{\bm{k}}Q_{\bm{k}}(\bm{x}) \label{psiG}\\
  && d\mu_{\psi} = \psi^2(\bm{x})\omega(\bm{x})d\bm{x} \label{dPgmin}\\
  &&\frac{\Braket{g(\bm{x})\psi^2(\bm{x})}}{\Braket{\psi^2(\bm{x})}} \to \min
  \label{Pgmin} \\
  &&\frac{\sum_{{\bm{jk}}}\alpha_{\bm{j}} \Braket{g(\bm{x})Q_{\bm{j}}(\bm{x})Q_{\bm{k}}(\bm{x})}\alpha_{\bm{k}}}{\sum_{{\bm{jk}}}\alpha_{\bm{j}}\Braket{Q_{\bm{j}}(\bm{x})Q_{\bm{k}}(\bm{x})}\alpha_{\bm{k}}} \to \min
  \label{PgminQ}
\end{eqnarray}
Instead of solving (\ref{gmin}): find $\bm{x}$ providing minimal $g$,
solve (\ref{Pgmin}) instead: find  probability state (a wavefunction $\psi(\bm{x})$) (\ref{psiG}),
providing minimal expected $g$ (\ref{Pgmin}).
After  expansion to (\ref{PgminQ}) obtain 
generalized eigenvalues problem (\ref{GEV}) 
with $M^{L}_{\bm{jk}}=\Braket{g(\bm{x})Q_{\bm{j}}(\bm{x})Q_{\bm{k}}(\bm{x})}$ and
$M^{R}_{\bm{jk}}=\Braket{Q_{\bm{j}}(\bm{x})Q_{\bm{k}}(\bm{x})}$,
that can be efficiently solved numerically\cite{lapack}.
The result is $\psi^{[\min]}(\bm{x})$,
corresponding to minimal eigenvalue $\lambda^{[\min]}$,  providing
probability distribution (\ref{dPgmin}),
not some specific $\bm{x}$ value as when solving  (\ref{gmin}) problem directly.
The answer in a form of probability distribution
is typically the most convenient in applications (because it allows to decouple observable $g$ (or $\bm{x}$)
and probability)
and used in other techniques, such as Bayesian Learning.
If, for any reason,
the $\bm{x}$, corresponding to found probability distribution
is required, it can be estimated as
\begin{eqnarray}
  \bm{x}^{[\min]}_{est}&=&\frac{\Braket{\left(\psi^{[\min]}(\bm{x})\right)^2\bm{x}}}{\Braket{\left(\psi^{[\min]}(\bm{x})\right)^2}}
  \label{xest}
\end{eqnarray}
and global minimum of  $g$ can be estimated as $\lambda^{[\min]}$.
The  $\bm{x}^{[\min]}_{est}$, besides being interpreted as optimization answer,
can also be used by other optimization algorithms as starting value.
Another interesting research topic is the roots
of  $\psi^{[\min]}(\bm{x})$. The value of $g$ is large
near $\psi^{[\min]}(\bm{x})$ roots, which typically
correspond to the ``spikes'' in $g$.
This is especially simple in 1D case
with polynomial basis\cite{2015arXiv151005510G}: for a given $n$
the $\psi^{[\min]}(x)$ ($n-1$ order polynomial) has exactly $n-1$
simple real distinct roots (but not necessary on the support
of the measure $d\mu=\omega(x)dx$).

Other than $\psi^{[\min]}(\bm{x})$ eigenvectors of (\ref{GEV}) are also
of interest in applications. First order variation (\ref{var1G})
of Rayleigh quotient 
is equal to zero, because of (\ref{GEV}).
Second order variation of Rayleigh quotient has a simple form  (\ref{var2G}),
which is positive for arbitray $\delta\psi$ when $i=\min$ (global minimum $\lambda^{[\min]}$).
\begin{eqnarray}
  \psi&=&\psi^{[i]}(\bm{x})+\delta\psi 
  \label{varpsi} \\
\frac{\Braket{g(\bm{x})\psi^2}}{\Braket{\psi^2}}
&=& \lambda^{[i]} \label{var0} \\
&+&  2\left[\Braket{g(\bm{x})\psi^{[i]}\delta\psi} -
    \lambda^{[i]}\Braket{\psi^{[i]}\delta\psi}
      \right]      
      \label{var1G}\\
      &+&\left[\Braket{g(\bm{x})\left(\delta\psi\right)^2} -
         \lambda^{[i]}\Braket{\left(\delta\psi\right)^2}\right] + \dots
     \label{var2G} 
\end{eqnarray}
All found $\psi^{[i]}(\bm{x})$ states (\ref{psiC})
correspond to $\min$, $\max$ or saddle point of $g$
(first order variation (\ref{var1G}) is zero)
with respect to wavefunction variation (\ref{varpsi}).
Thus this technique can be considered as a type
of differential calculus (generalized Radon--Nikodym).
This calculus variate probability state (wavefunction (\ref{varpsi})),
not $g(\bm{x}+\delta \bm{x})$ argument,
as regular calculus does. This transition
from  $\bm{x}+\delta \bm{x}$ variation to wavefunction
 variation (\ref{varpsi})
allows to transform
original optimization problem (\ref{gmin}) in $\bm{x}$ space
to Rayleigh quotient optimization (\ref{PgminQ}) in $\psi(\bm{x})$
space and then to
generalized eigenvalues
problem (\ref{GEV}), that can be efficiently solved numerically\cite{lapack}.

The approach of this paper works with (\ref{dPgmin}) probability distribution.
Bayesian Learning\cite{theodoridis2015machine}
works with $\rho(\bm{x}|\theta)$ parametric distribution.
Despite both approaches work with probability distributions in $\bm{x}$ space,
there are important conceptual differences between (\ref{Pgmin}) approach and Bayesian Learning:
\begin{itemize}
\item Probability densities are represented \textsl{only in} (\ref{dPgmin}) form.
  Original argument $\bm{x}$ enter to optimization problem (\ref{Pgmin})
  only via $\psi(\bm{x})$ components (\ref{psiG}).
\item There is no any norm or standard deviation $\sigma^2$ involved,
  what make our approach applicable to strongly non--Gaussian case.
\item Unique ``natural basis'', the eigenvectors of (\ref{GEV}) problem,
  can be considered as cluster centers.
  Predicted value can now be obtained as a superposition of projections\cite{malyshkin2015norm}
  of these cluster centers to the state of interest, such as in Eq. (\ref{RNsimpleEVBasis})
  for interpolation problem.
\item The $\bm{x}^{[\min]}_{est}$ answer (\ref{xest})
  obtained using (\ref{dPgmin}) distribution with  $\psi^{[\min]}$
  solution of (\ref{GEV}),
  is some kind similar to $\rho(\bm{x}|\theta)$ parametric distribution used in
  Bayesian Learning\cite{theodoridis2015machine},
  but is obtained from very different considerations.  
\item The (\ref{GEV}) give global optimization solution $\psi^{[\min]}$
  and the minimum is $\lambda^{[\min]}$.
\end{itemize}

Developed software (see Appendix \ref{appnum} for description),
besides the main goal of this study, obtaining relaxation rates distribution, 
can be also used for 1D optimization using the technique described.
Output \verb+*_spectrum.dat+ files
contain three columns: (index $i$, $\lambda^{[i]}$, $x^{[i]}_{est}$),
calculated for all eigenstates $i=[0\dots n-1]$ of (\ref{GEV}).
The minimal eigenvalue ($i=0$) corresponds to the value of global minimum
and corresponding $x^{[i]}_{est}$ is the (\ref{xest}) estimation.

There is an important extension of  original problem
from using input data in a form of timeserie (\ref{tsdata}),
with $\bm{x}$ and $f$ being  
experimentally measured input variables,
to 
the theory of \textsl{distribution regression problem\cite{dietterich1997solving}},
where a bag of $\bm{x}$ observations, not a single $\bm{x}$ realization,
is mapped to a single $f$:
\begin{eqnarray}
  (\bm{x}_1,\bm{x}_2,\dots,\bm{x}_j,\dots,\bm{x}_N)_{l}&\to&f_l  \label{regressionproblem} \\
  l&=& 1\dots M \nonumber
\end{eqnarray}
The (\ref{regressionproblem}) can be considered as a distribution of $\bm{x}$
(a bag with $N$ realizations of $\bm{x}$ used as input)
is  mapped to a single $f$ value.
The theory of above
can be transformed\cite{2015arXiv151107085G,2015arXiv151109058G} to use a distributions of $\bm{x}$ as input.
In \cite{2015arXiv151109058G}
a generalization of formulas
(\ref{RNsimple}), (\ref{RNsimpleEVBasis}) and (\ref{GEV})
is obtained for this case.
For distribution regression problem
the (\ref{dPgmin}) have a meaning  of ``the distribution
of distribution'' and can be effectively used
in applications,
such as
to study the data with noise,
the examples of  Ref. \cite{yang2005review}
(drug activity prediction,
content--based image classification,
text categorization),
and studied in Ref. \cite{taleb2014silent}
the concept of ``volatility of volatility''
that directly correspond to the distribution
of distribution solution\cite{2015arXiv151109058G}
of distribution regression problem.
Similary,
in case of global optimization in distribution regression problem,
the solution of (\ref{Pgmin})
give \textsl{the distribution of $\bm{x}$--distribution},
(\ref{dPgmin}), providing minimal $g$.

\section{\label{RNSpectdem} Demonstration Of Radon--Nikodym Spectral Approach}
The main idea of Radon--Nikodym Spectral Approach
is to construct, from experimental observations,
a ``Quantum Mechanics Hamiltonian'', the spectrum of which
correspond to the dynamic characteristic of interest,
 then to study this Hamiltonian spectrum
considering obtained eigenvalues distribution
as a ``density of states''.
Fourier analysis is different,
it works with
interpolation of an observable value (using $L^2$ norm)
and harmonics
are initially selected, not obtained from the data.
Laplace analysis, in addition to the problems of Fourier analysis,
has the problems of sample insufficient size and discretization noise,
thus is not applicable to sampled data.

Radon--Nikodym Spectral Approach,
same as quantum mechanics,
constructs probability states,
then observable variable characteristics
are calculated using obtained probability densities.
There is no any $L^2$ norm involved, what allows to study
even strongly non--Gaussian distributions\cite{malyshkin2015norm}.
The main idea is to obtain from experimental data, such as
(\ref{Cqm}) and  (\ref{CDqm}),
two {\bf matrices} (not {\bf vectors} like Fourier components
or least squares coefficients).
Then solve generalized eigenvalues
problem (\ref{GEV}) with these two matrices
(diagonalize them simultaneously).
Obtained eigenvalues spectrum  provide
the distribution of an observable.
In random matrix theory\cite{guhr1998random}
the goal is similar: to obtain
matrix eigenvalues density,
but the matrix in study is typically obtained from some initially selected model\cite{erdos2016dynamical}.
In this paper the $M^{L}_{jk}$ and $M^{R}_{jk}$ matrices, determining virtual Hamiltonian,
are obtained directly from sampled signal data.
 Similary to random matrix theory
eigenvalues density determine the
distribution of the characteristic in study.
Few  examples
of $M^{L}_{jk}$ and $M^{R}_{jk}$ choice
(for practical implementation only one matrix
$M^{L}_{jk}$ or $M^{R}_{jk}$ need to be positively defined).
\begin{itemize}
    \item
If $\Braket{Q_j df/dx Q_k}$ is used as $M^{L}_{jk}$
and $\Braket{Q_j Q_k}$ is used as $M^{R}_{jk}$ in (\ref{GEV})
then the distribution of obtained $\lambda^{[i]}$; $i=[0..n-1]$ 
is the distribution of  $df/dx$, observed in the data.

  \item
If $\Braket{Q_j df/dx Q_k}$ is used as $M^{L}_{jk}$
and $\Braket{Q_j f Q_k}$ is used as $M^{R}_{jk}$ in (\ref{GEV})
then the distribution of obtained $\lambda^{[i]}$; $i=[0..n-1]$ 
is the distribution of relaxation rates $\frac{df}{f dx}$
observed in the data.

  \item
If $\Braket{Q_j \frac{d\ln(f)}{dx} Q_k}$ is used as $M^{L}_{jk}$
and $\Braket{Q_j Q_k}$ is used as $M^{R}_{jk}$ in (\ref{GEV})
then the distribution of obtained $\lambda^{[i]}$; $i=[0..n-1]$ 
is the distribution of relaxation rates $\frac{df}{f dx}$
observed in the data.

\item
If $\Braket{Q_j x Q_k}$ is used as $M^{L}_{jk}$
and $\Braket{Q_j Q_k}$ is used as $M^{R}_{jk}$ in (\ref{GEV}),
then obtained eigenvalues $\lambda^{[i]}$; $i=[0..n-1]$ 
are the nodes
and $1/\left(\psi^{[i]}(\lambda^{[i]})\right)^2$
are the weights
of $n$--point Gauss quadrature
built on $d\mu=\omega(x)dx$ measure\cite{2015arXiv151005510G}
(also note that $\psi^{[i]}(x)$ from (\ref{psiC})
are proportional to Lagrange interpolating polynomials:
$\psi^{[i]}(\lambda^{[j]})=0$ for $i\ne j$).
Because of this fact the software from appendix \ref{appnum}
can be used to build Gauss quadratures from sampled data.
\end{itemize}
Below we are going to demonstrate this spectral approach on several
examples for both model and real life data.
In  practice relaxation data is most often obtained
from the processes of  degradation, relaxation, market dynamics, etc.
We chose the following data for demonstration: 
\begin{itemize}
\item Li-Ion degradation, subsection \ref{TwoStage}.
  This data is relatively easy to obtain experimentally,
  and, using accelerated cycling, degradation signals can be
 recorded in a period of few months.
  This make the system a good testbed degradation example,
  compared to other systems where degradation
  processes take decades.
\item Supercapacitor discharge, subsection \ref{MultiStageSuperC}.
  Relaxation data in this system can be obtained in
  minutes, and, very important, the process is repeatable.
  In contrast with degradation process
  charge/discharge cycle can be repeated
  many times.
\item Stock price changes, subsection \ref{StockP}.
  It is well known\cite{mandelbrot2014misbehavior}
  that financial markets do not have any short--term
  time--scale and the distribution of market scales
  is actively studied\cite{taleb2014silent}.
  This make interesting to apply our theory to
  market data, because we can obtain the distribution of price change rates
  directly from timeserie sample.
\end{itemize}

\subsection{\label{TwoStage} Li-Ion Degradation}

\begin{figure}[!t]
  \includegraphics[width=4cm]{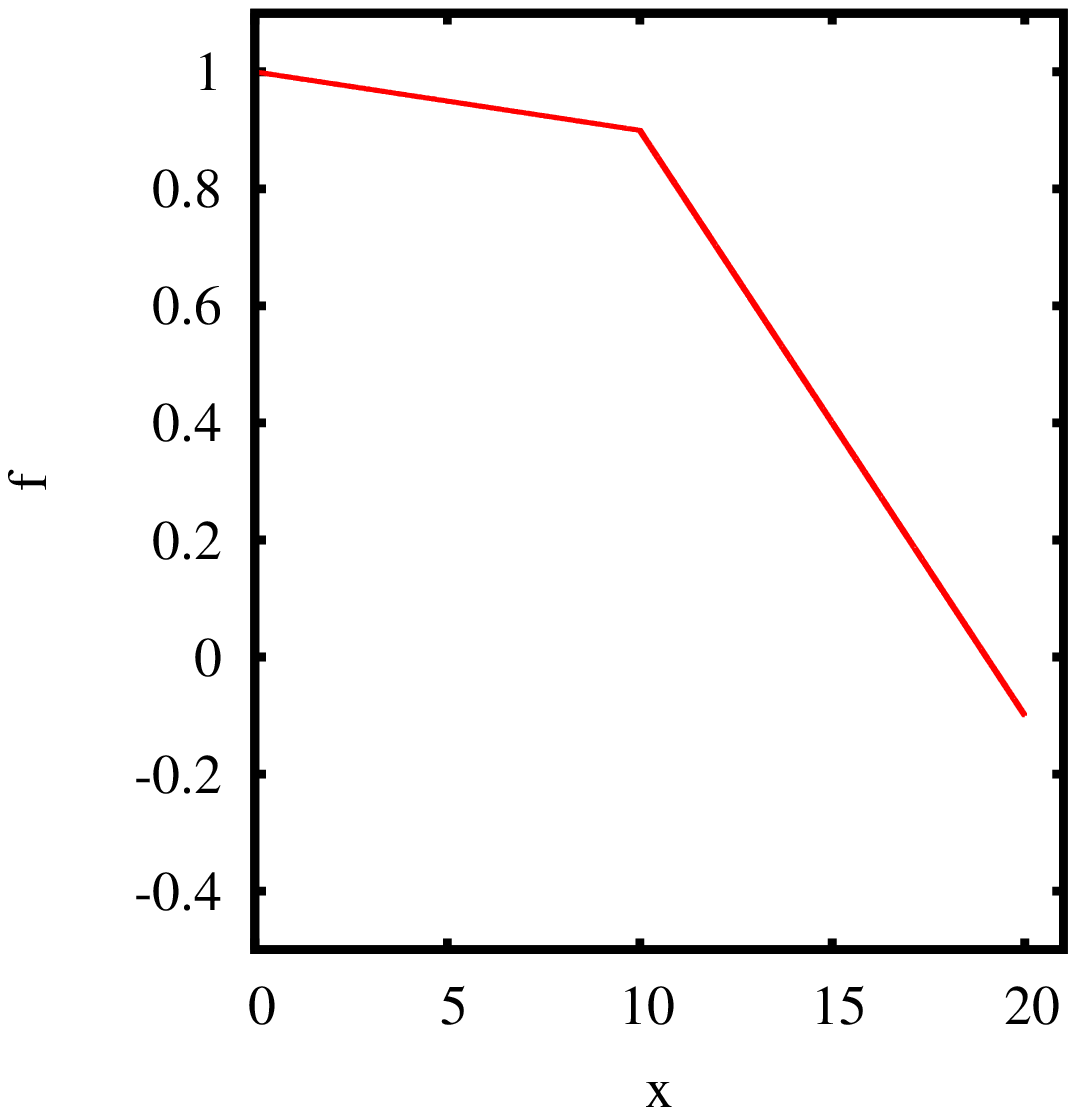}
  \includegraphics[width=4cm]{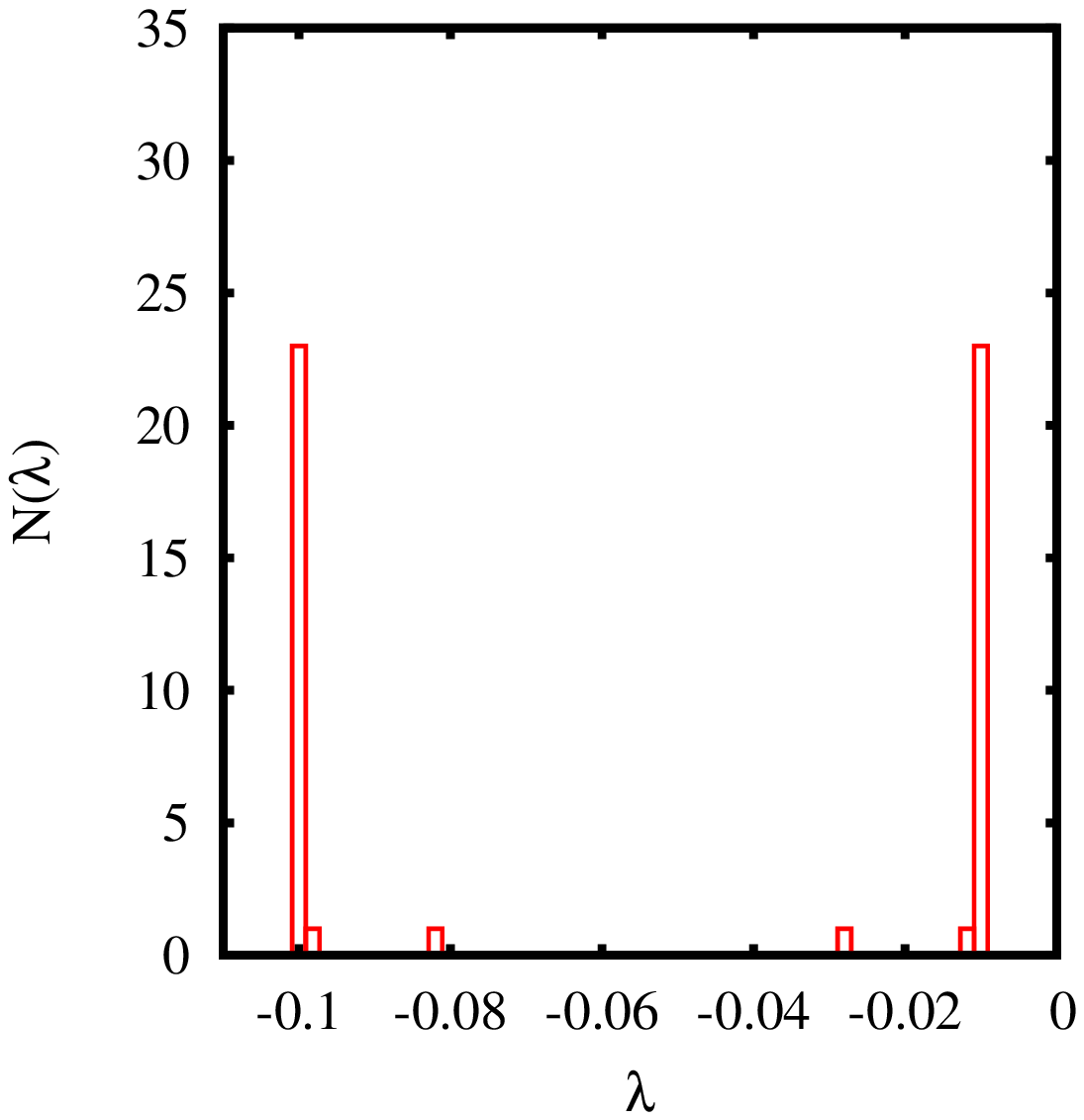} \\
  \includegraphics[width=4cm]{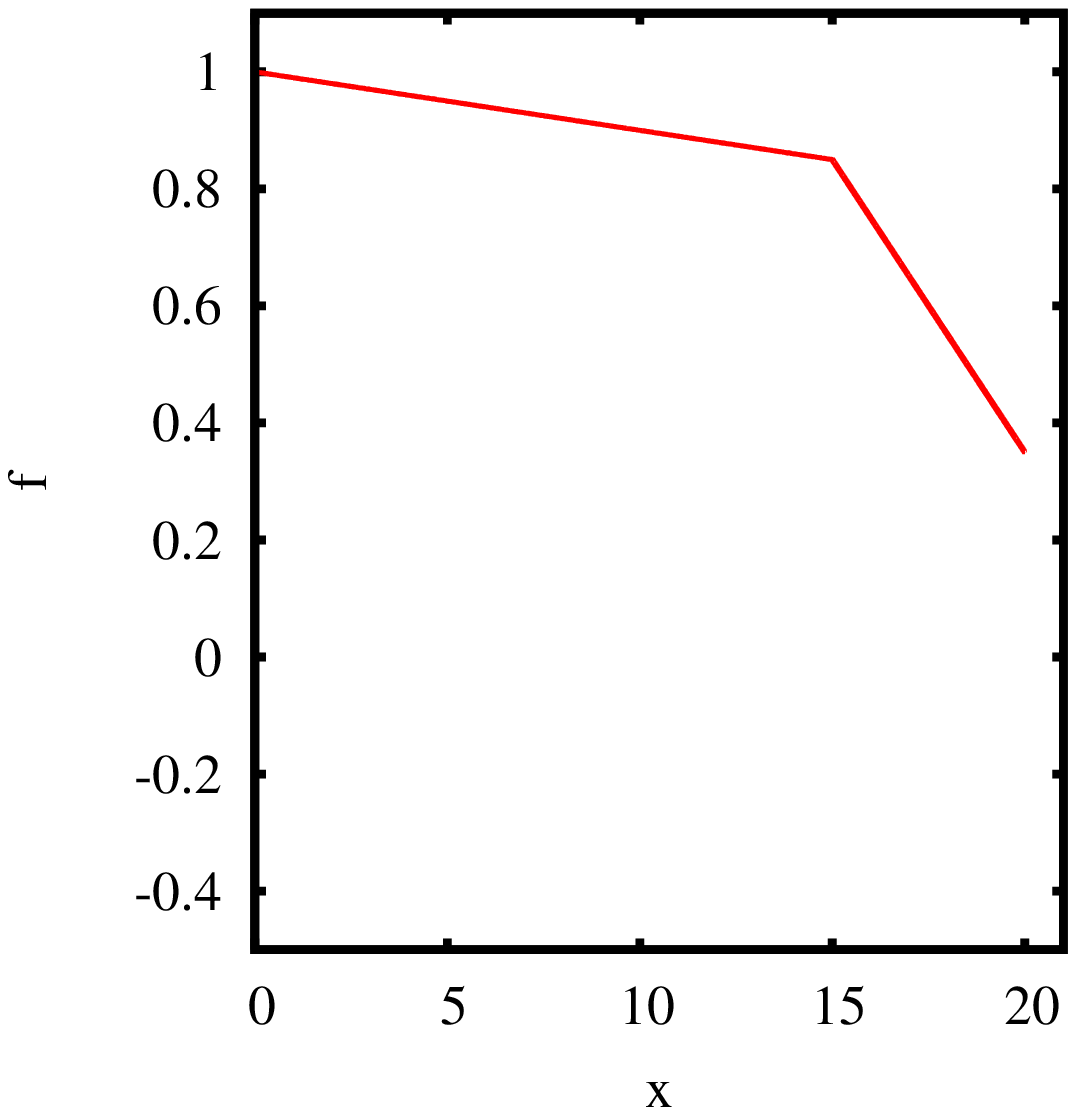}
  \includegraphics[width=4cm]{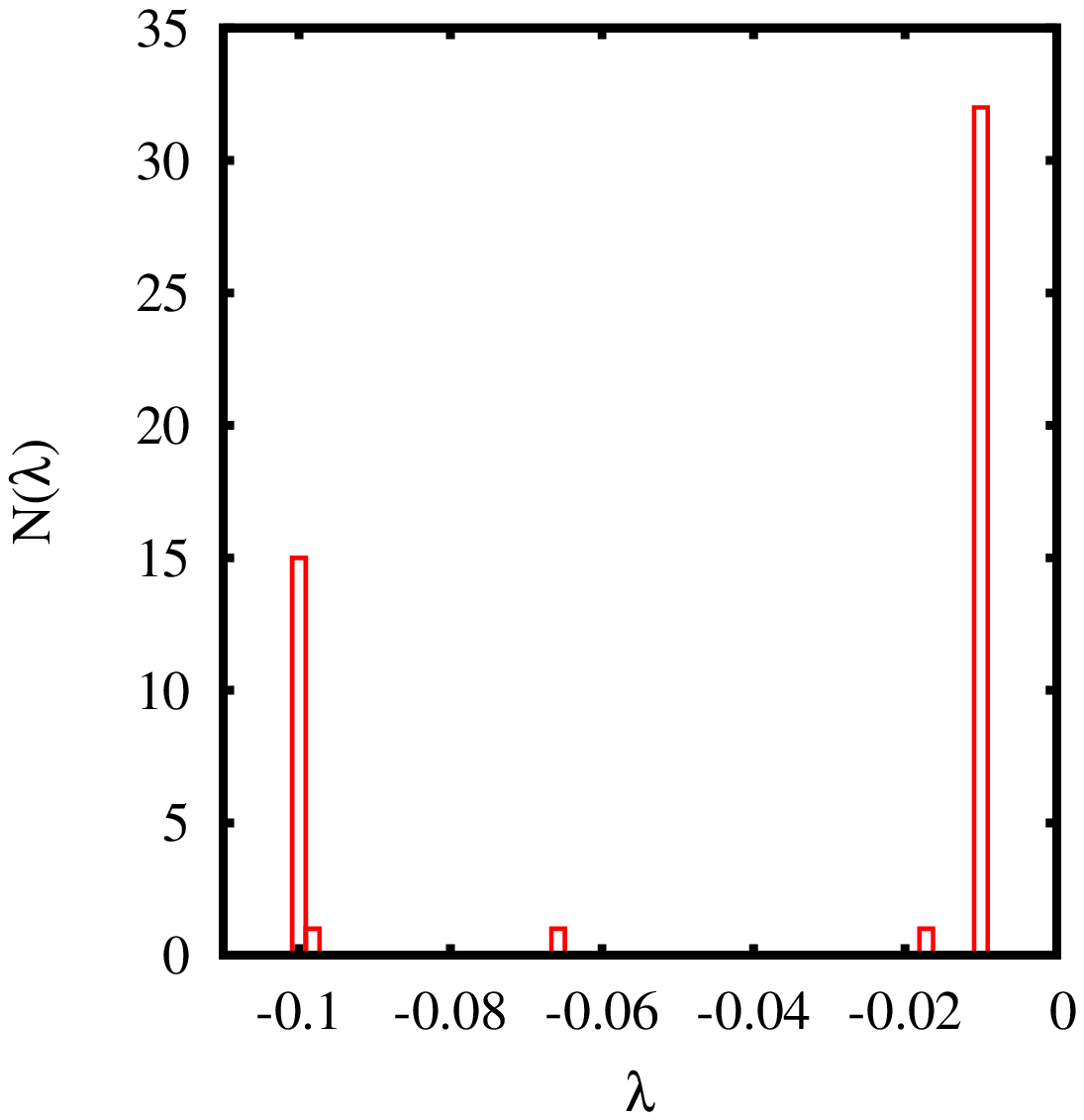}
\caption{\label{figLiIon10}
  Two stage degradation Li-Ion model with
  the slope on first and second stages $-0.01$ and $-0.1$ respectively.
  The stages length is 10:10 (equal time) for top chart
  and 15:5 for bottom chart.
  Corresponding distributions of $\lambda$ are calculated with $n=50$.
  }
\end{figure}

\begin{figure}
  \includegraphics[width=4cm]{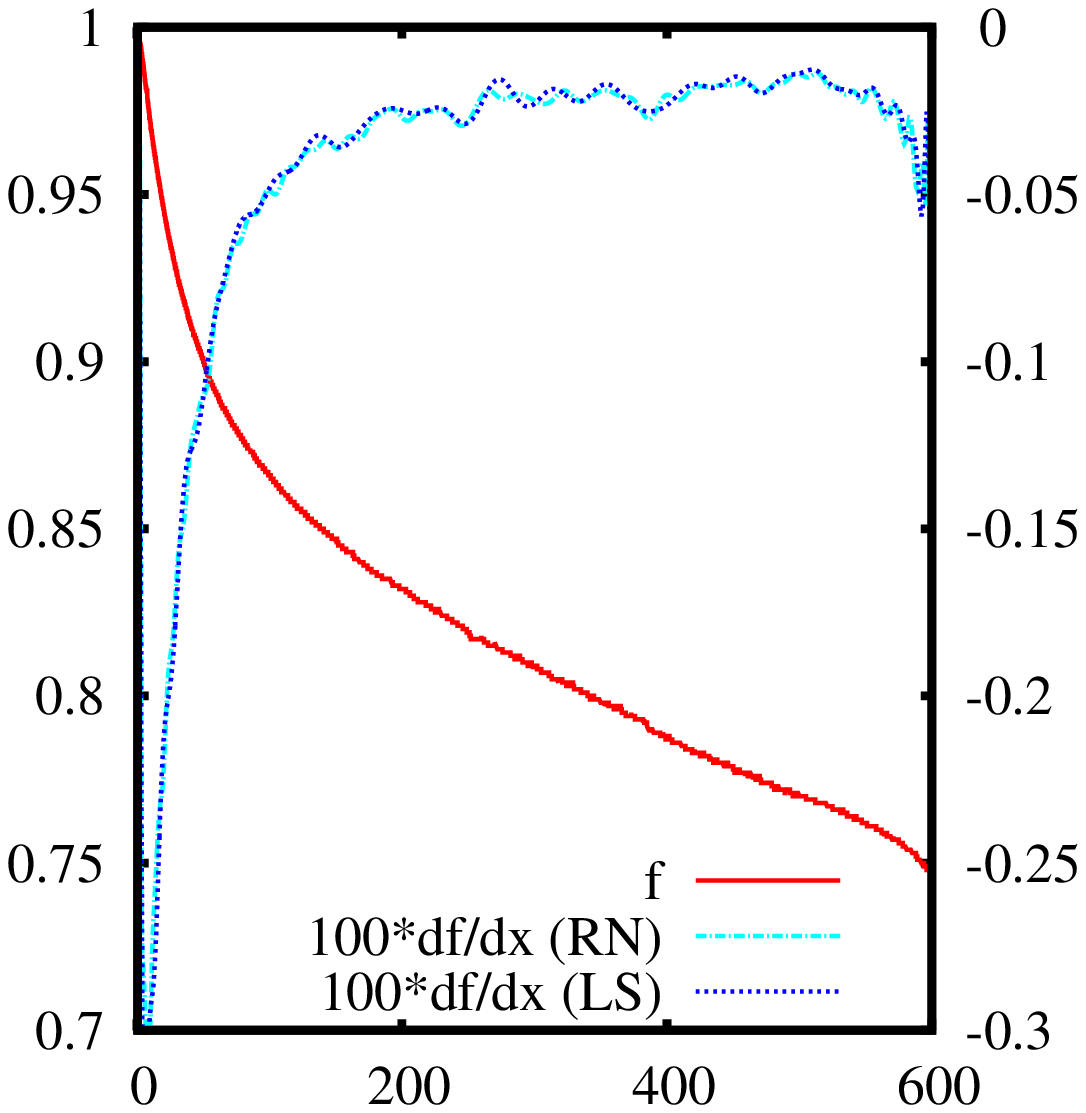}
  \includegraphics[width=4cm]{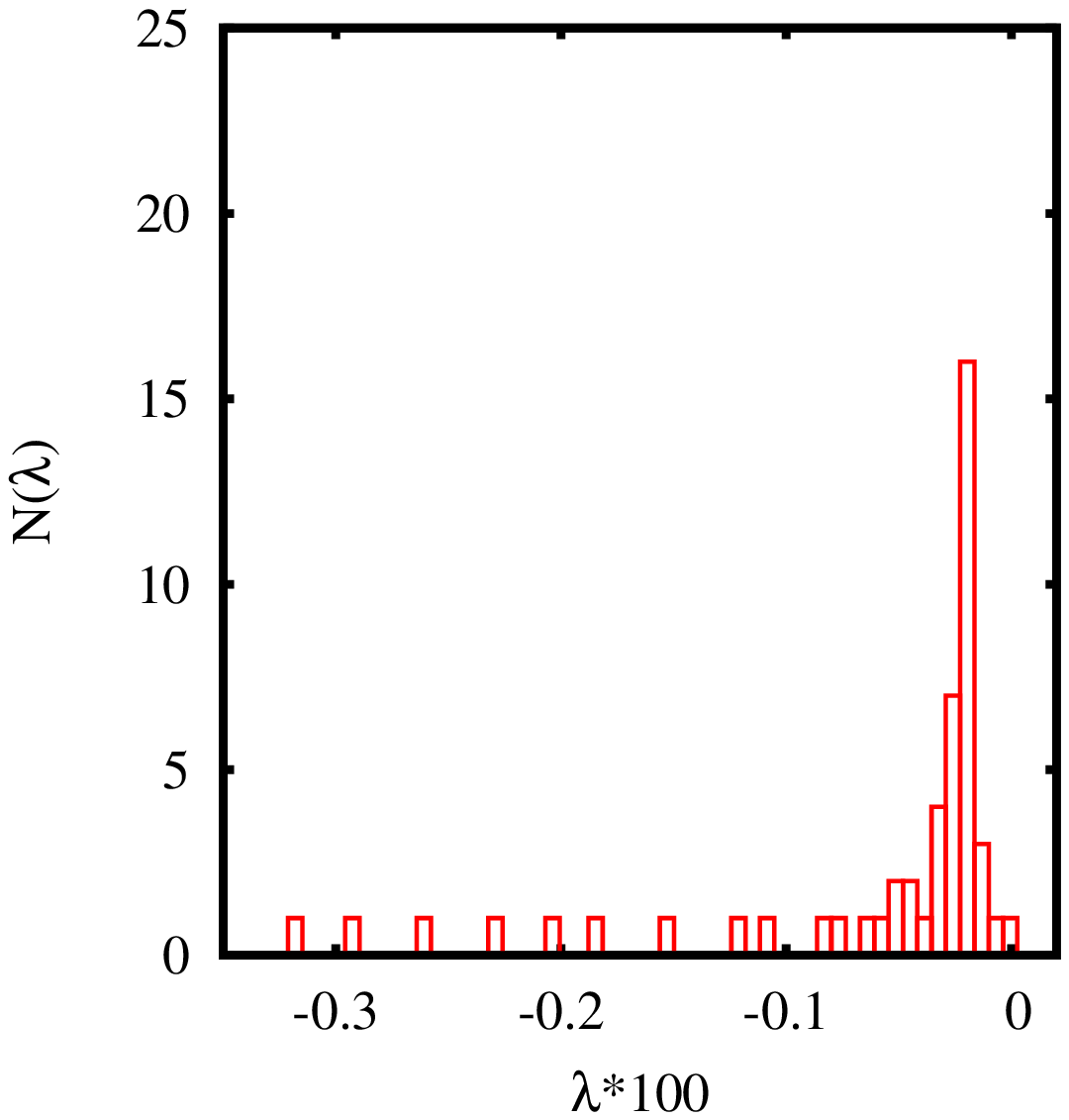}
\caption{\label{figLiIonLG}
  Degradation of LG Chem\cite{LGChemLGHG2} Li-Ion battery.
  Left: dependence of battery capacity $f$ on cycle number $x$,
  and interpolation of degradation rate $df/dx$ using least squares (\ref{leastsq})
  and Radon--Nikodym (\ref{RNsimple}) interpolation
  for $df/dx$.
  Corresponding distributions of $\lambda$ are calculated with $n=50$.
  }
\end{figure}

Consider typical for degradation dynamics a
simple two--stage Li-Ion degradation model\cite{MalLiIon2016,liionizversiyaran},
when battery capacity $C$ linearly decay with cycle number $N$,
but the slope (degradation rate) changes at some point, see Fig. \ref{figLiIon10}.
To build a ``Hamiltonian'' the spectrum of which
give degradation rate take $f=C$, $x=N$, $g=df/dx$,  
two matrices need to be calculated from experimental observations:
$\Braket{Q_j df/dx Q_k}$ (\ref{CDqm}) to use it as $M^{L}_{jk}$
and $\Braket{Q_j Q_k}$ (\ref{qm}) as $M^{R}_{jk}$ in (\ref{GEV}).
Then $\lambda$ has exactly the meaning of degradation rate $df/dx$,
and the distribution of it carry the information
about available degradation rates in experimental observations.
In Fig. \ref{figLiIon10} the distributions
of $\lambda$, assuming all the eigenvalues have equal weights,
are presented.
These are  two mode distributions
with peaks at exact degradation rate
of first and second stages: $\lambda=-0.01$ and $\lambda=-0.1$.
Peak height growths with the increase of observations number at
 degradation rate value, but the dependence
is more complicated than simple ratio of observations number.
The relation was later found in Ref. \cite{ArxivMalyshkinLebesgue}
where the concept of Lebesgue integral quadrature
was introduced. Each eignevalue $\lambda^{[i]}$
is considered as \textbf{value--node} $g_i$ of Lebesgue quadrature
and corresponding weight $w_i$ are:
\begin{eqnarray}
  \Braket{g}&=&\sum\limits_{i=0}^{n-1} g_i \Braket{\psi^{[i]}}^2 
  \label{inegralGsum} \\
  g_i&=&\lambda^{[i]} \label{fiLeb} \\
  w_i&=& \Braket{\psi^{[i]}}^2 \label{wiLeb}
\end{eqnarray}
If instead of equal weights for $\lambda^{[i]}$ as in Fig. \ref{figLiIon10}, one take Lebesgue quadrature wights (\ref{wiLeb})
then, because we have chosen $d\mu=dN$ measure,
distribution peaks height correspond exactly to
stage length\cite{ArxivMalyshkinLebesgue}.
Lebesgue quadrature (\ref{inegralGsum}
can be considered as Lebesgue integral interpolating formula,
by $n$--point discrete  measure,
the value--nodes $g_i$ select optimal positions of function values,
they are $\|g\|$ operator  eigenvalues (\ref{GEV}),
the weight $w_i$ is the measure
of all $g\approx g_i$ sets. Note, that weights (\ref{wiLeb}) give $\Braket{1}=\sum_{i=0}^{n-1}w_i$,
same as for Gaussian quadrature weights.

A more interesting data sample
is real world Li-Ion degradation data.
From \cite{LGChemLGHG2} datasheet,  battery capacity ($f$)
as a function of cycle number ($x$) is obtained.
Same as in the example above
we build $\Braket{Q_j df/dx Q_k}$ and
$\Braket{Q_j Q_k}$ matrices, from measured data,
then  use these matrices
in (\ref{GEV}) to obtain eigenvalues distribution.
Eigenvalues disribution with Lebesgue weights (\ref{wiLeb})
give spectrum distribution of ``what was already observed''
in degradation rates.
The distribution with equal weight
for each eigenvalue, Fig. \ref{figLiIonLG}, may serve
as an idicator of ``what can happen''
to degradation rates.
Strong difference in these distributions
possibly indicate future drastic changes in degradation rate.

Least squares (\ref{leastsq})
and Radon--Nikodym (\ref{RNsimple}) $df/dx$
interpolations are also calculated,
they are presented in Fig. \ref{figLiIonLG}.
The $df/dx$ data is smooth and both  least squares
and Radon--Nikodym  give very similar result in interpolation
of this specific example, not as with Lenna image in previous section.

\subsection{\label{MultiStageSuperC} Supercapacitor Discharge With 
  Relaxation Time Change}
\begin{figure}
  \includegraphics[width=4cm]{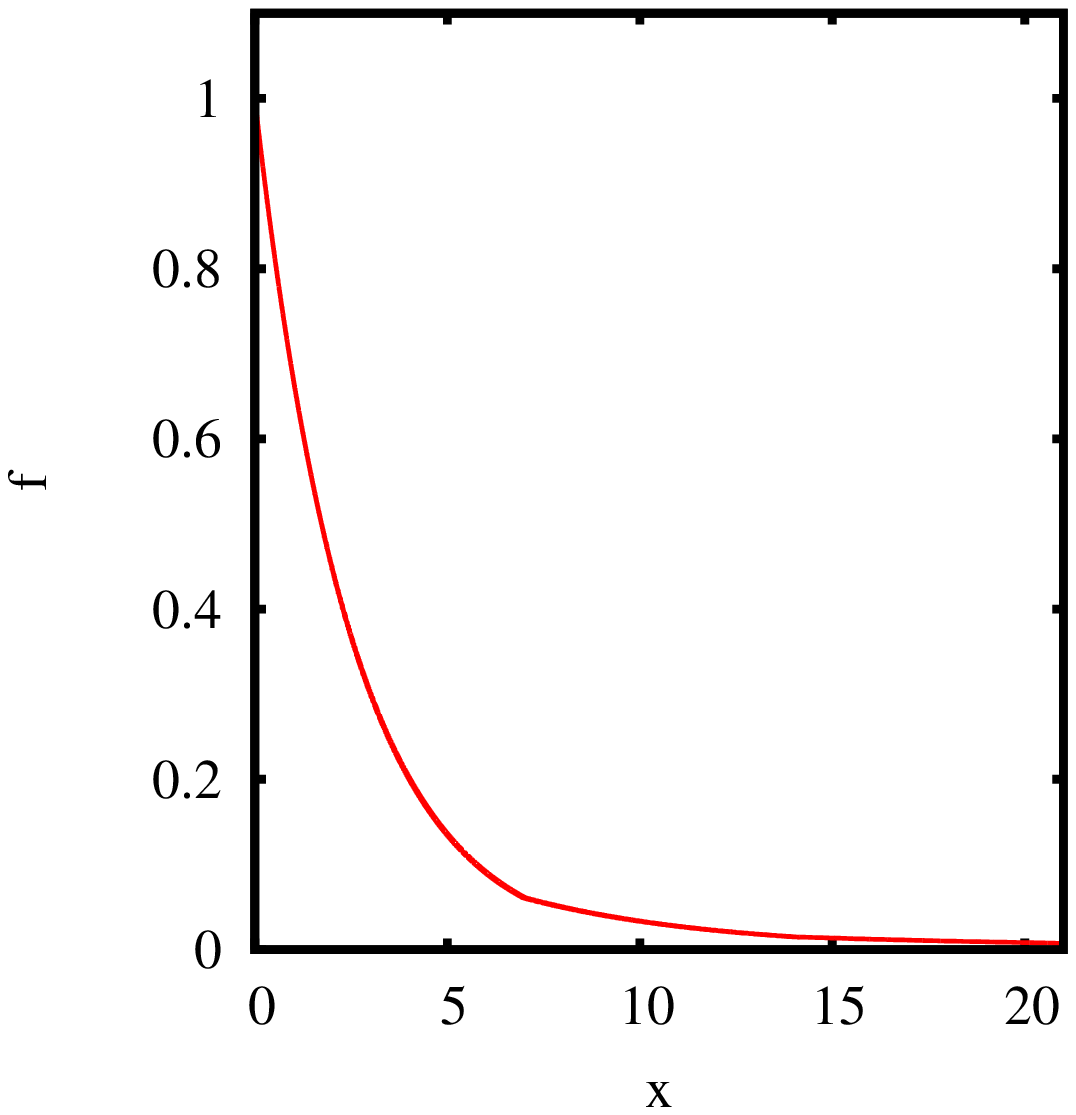}
  \includegraphics[width=4cm]{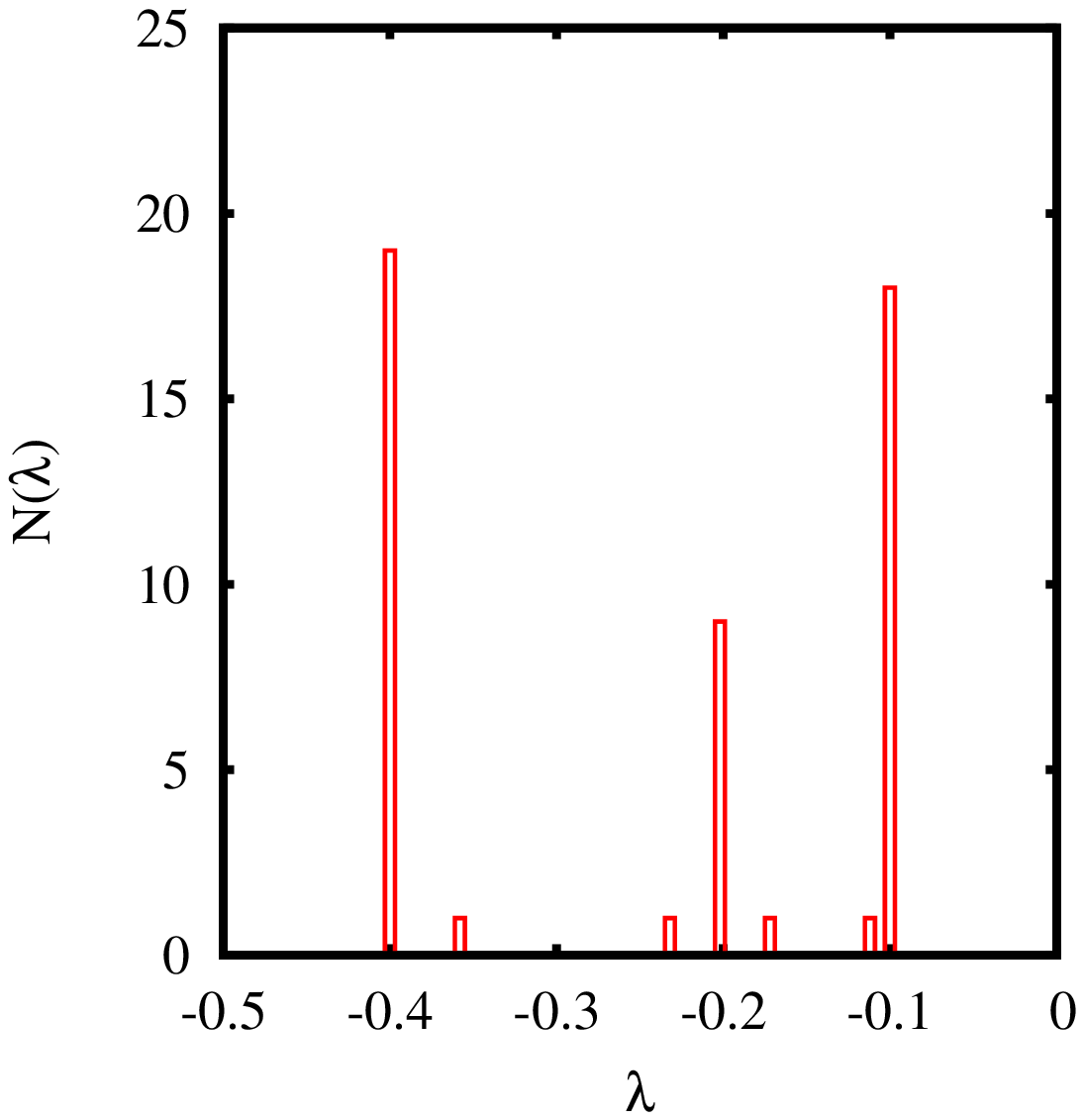} \\
  \includegraphics[width=4cm]{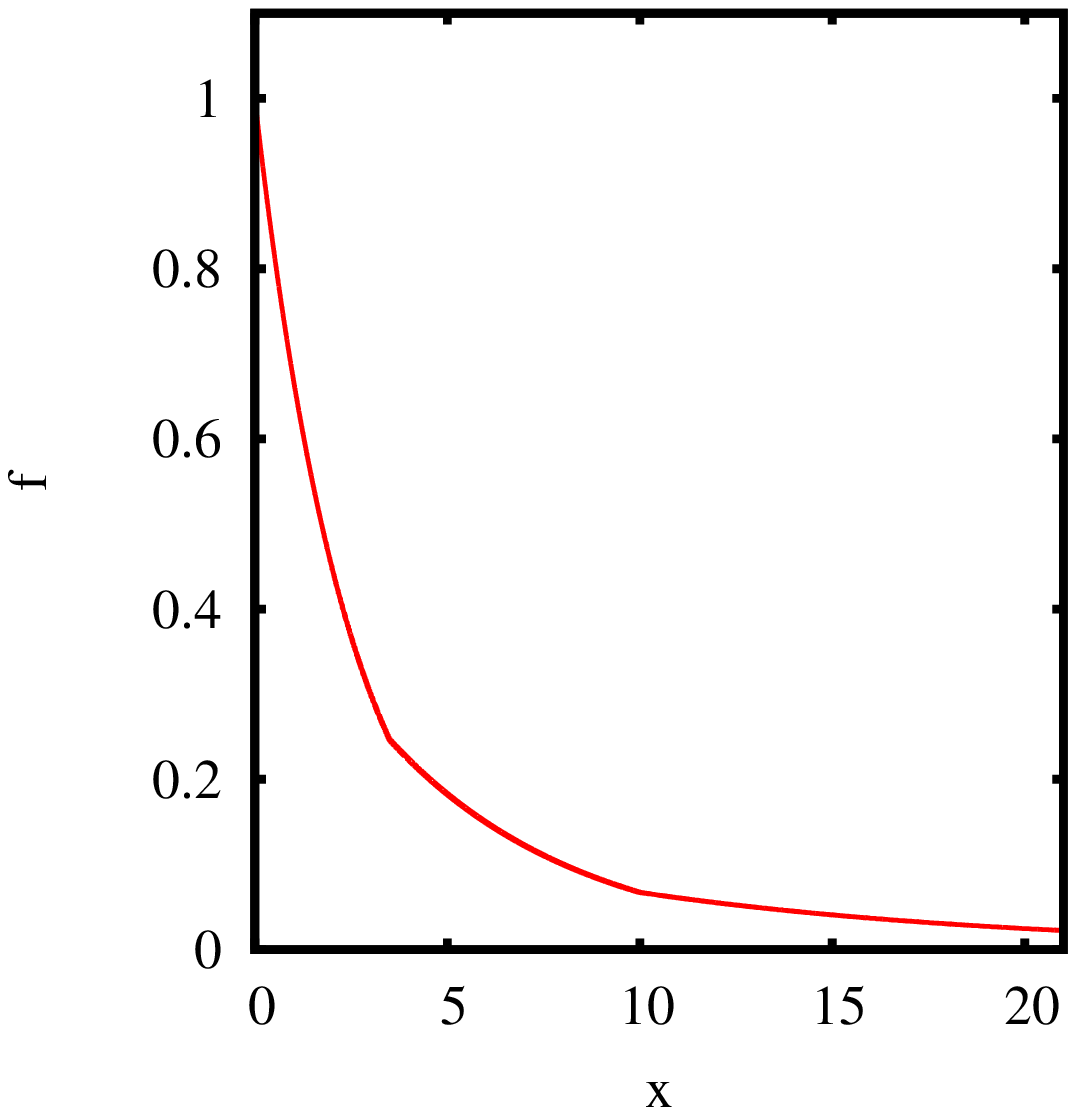}
  \includegraphics[width=4cm]{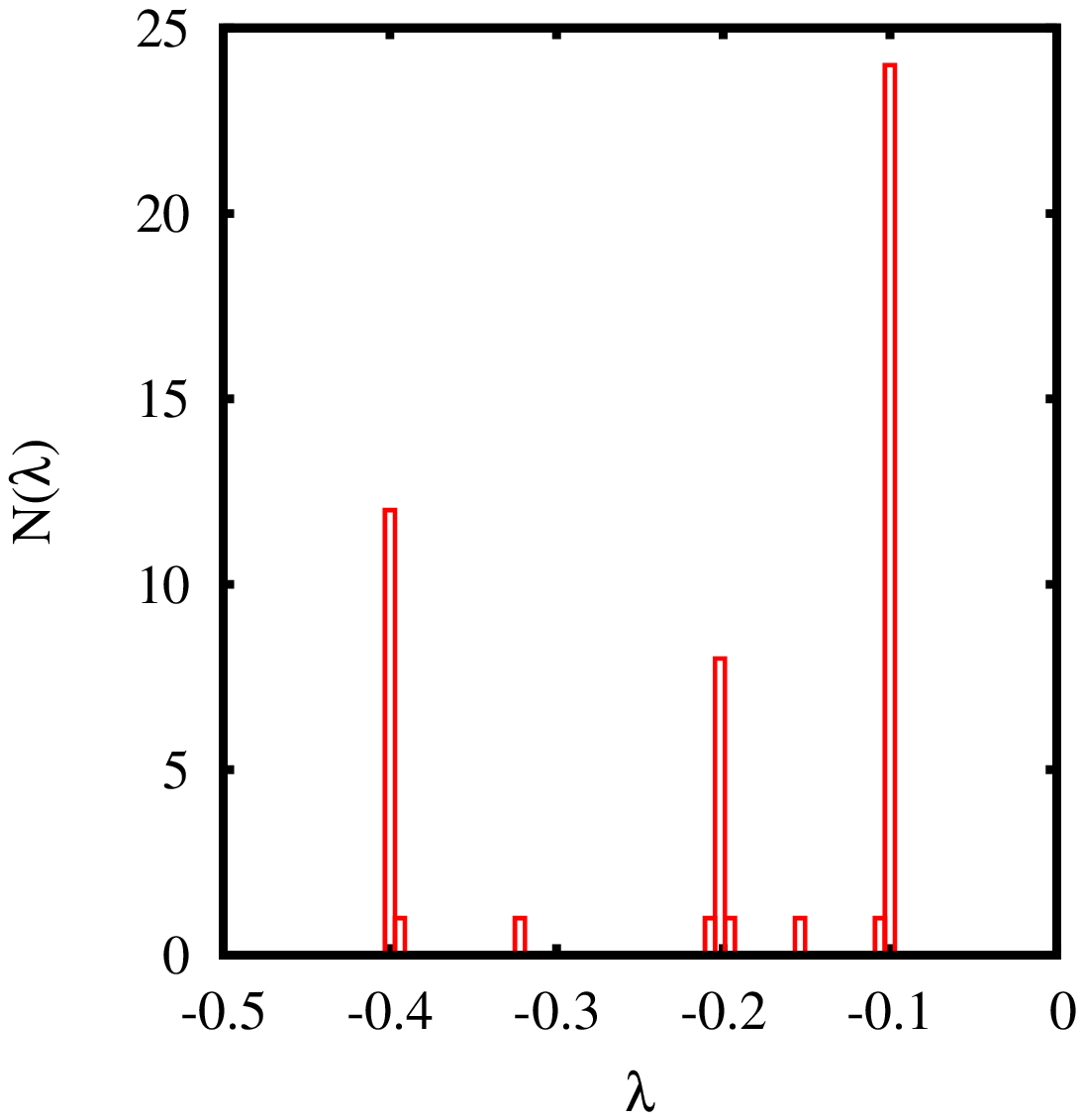}
\caption{\label{figSCDischarge}
  Three stage modeled supercapacitor discharge.
  Stages exponents  are $-0.4$, $-0.2$ and $-0.1$
  with (\ref{fpfL}) and (\ref{fpfR}) matrix selection.
  The stages length is 7:7:7 (equal time) for top chart
  and 3.5:7.0:10.5 for bottom chart.
  Corresponding distributions of $\lambda$ are calculated with $n=50$.
  }
\end{figure}

Consider a typical for relaxation
dynamics multi--stage supercapacitor discharge,
modeled with three consequential  relaxation times
as in Fig. \ref{figSCDischarge}.
There are two alternative 
of discharge rate calculation. First, corresponding to $f^{\prime}/f$, approach
\begin{eqnarray}
  M^{L}_{jk}&=& \Braket{Q_j \frac{df}{dx} Q_k} \label{fpfL}\\
  M^{R}_{jk} &=& \Braket{Q_j f Q_k} \label{fpfR}
\end{eqnarray}
and second, corresponding to $\ln(f)^{\prime}$, approach
\begin{eqnarray}
  M^{L}_{jk}&=& \Braket{Q_j \frac{d\ln(f)}{dx} Q_k} \label{lfL} \\
  M^{R}_{jk} &=& \Braket{Q_j Q_k} \label{lfR}
\end{eqnarray}
Both selections give $\lambda$ as a discharge rate.
The $\lambda$ has exactly the meaning of exponent ($\frac{df}{fdt}$)
and the distribution of it carry the information
about available exponents in experimental observations.

Consider first choice of  ``Hamiltonian''.
The matrices (\ref{fpfL}) and (\ref{fpfR}) are calculated from experimental observations as (\ref{CDqm}) and (\ref{Cqm}) to use them in (\ref{GEV});
In Fig. \ref{figSCDischarge} the distributions
of $\lambda$ are presented at right.
These are  three mode distributions
with peaks at exact exponents
on first, second and third stages: $\lambda=-0.4$,
$\lambda=-0.2$ and $\lambda=-0.1$.
Peak height growths with observations number at
the discharge rate, but, again, the dependence
is more complicated than simple ratio of observations number,
it is more similar to Christoffel function behavior\cite{malyshkin2015norm}.

\begin{figure}[h]
  \includegraphics[width=4cm]{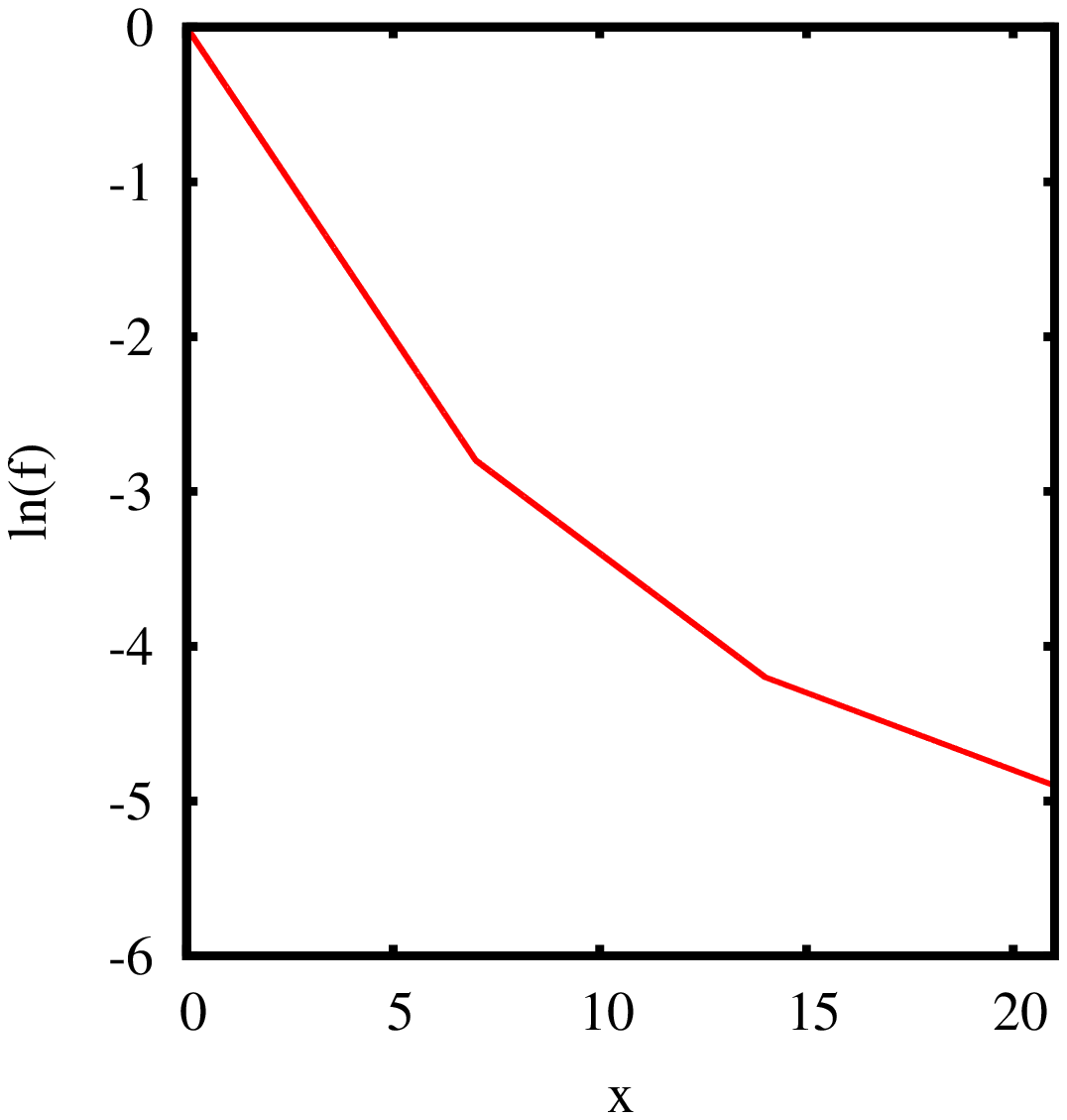}
  \includegraphics[width=4cm]{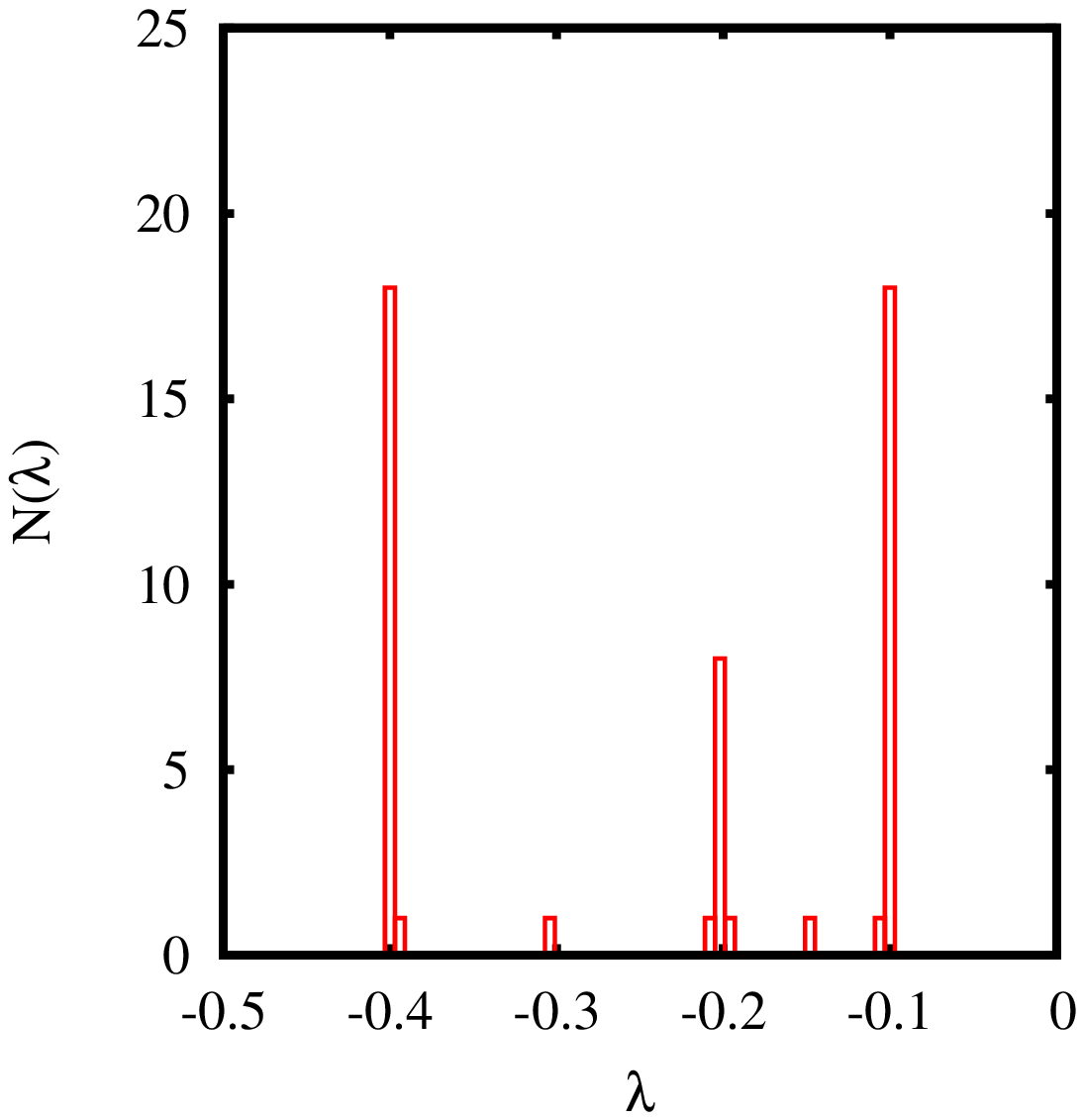} \\
  \includegraphics[width=4cm]{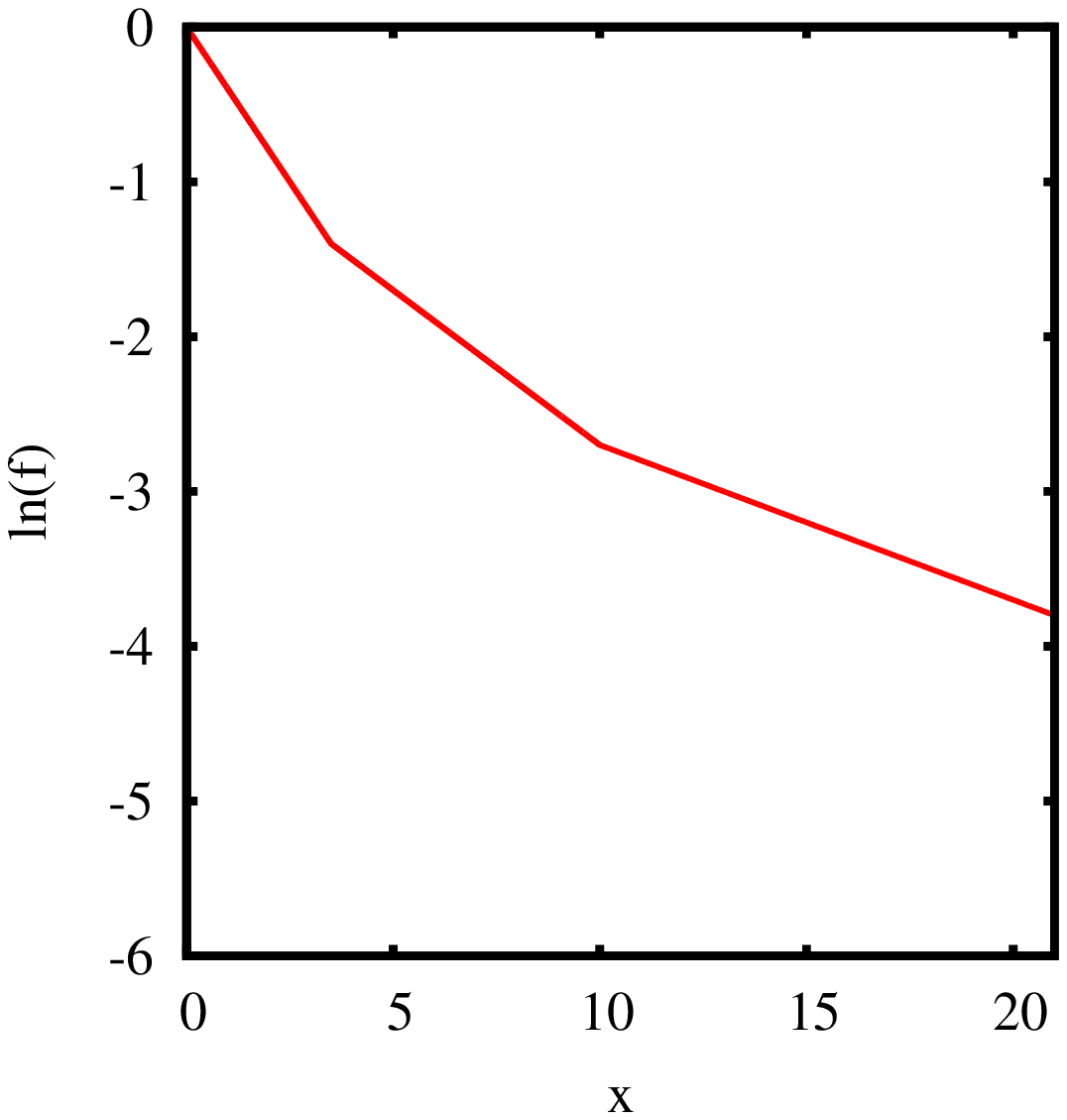}
  \includegraphics[width=4cm]{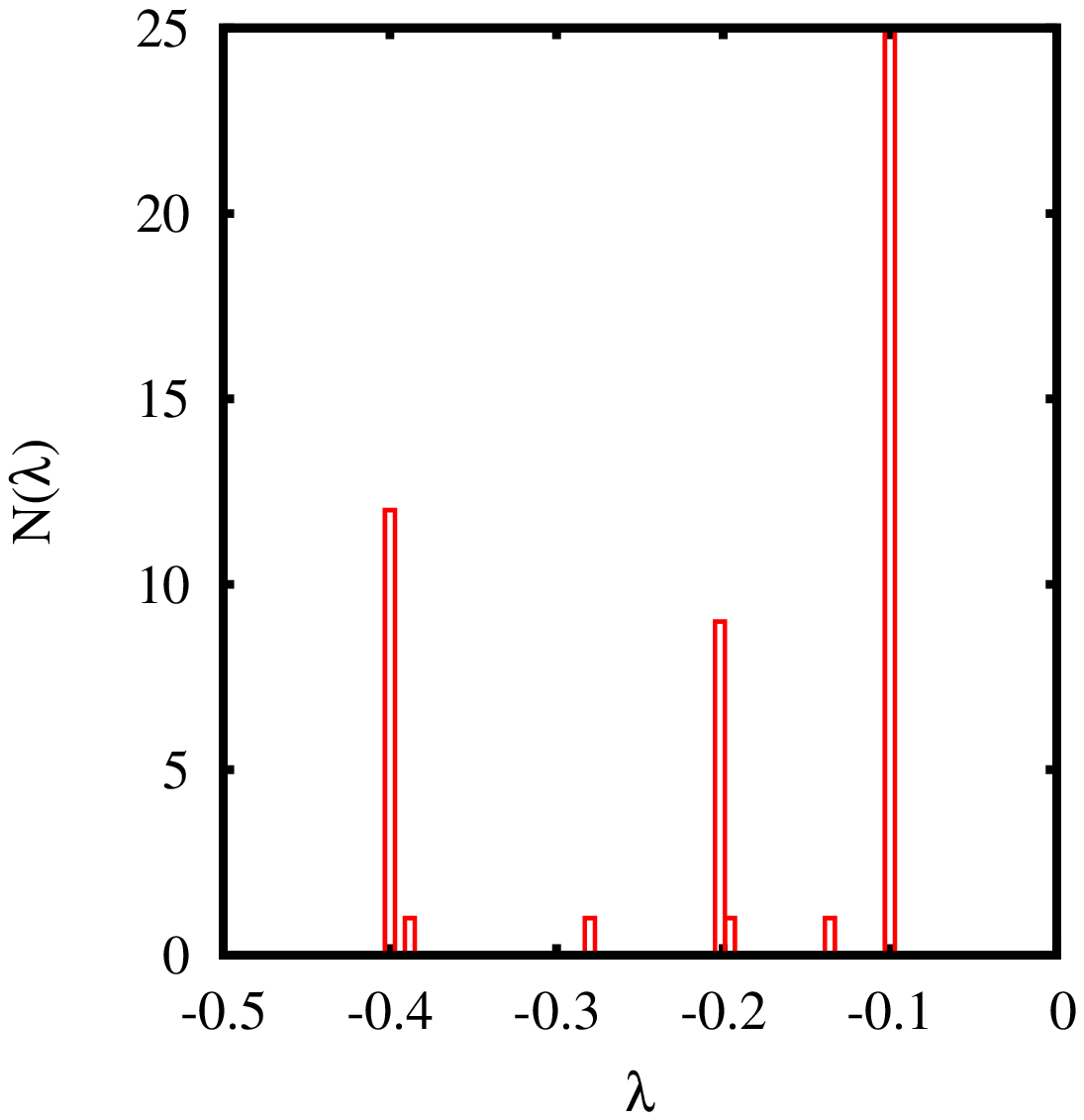}
\caption{\label{figSCDischargeLog}
  Three stage modeled supercapacitor discharge,
  same as in Fig. \ref{figSCDischargeLog},
  but with (\ref{lfL}) and (\ref{lfR}) matrix selection.
  }
\end{figure}

In Fig. \ref{figSCDischargeLog}
the results corresponding to second matrix choice 
(\ref{lfL}) and (\ref{lfR}) are presented, same approach as in section \ref{TwoStage} above.
The results in Figs. \ref{figSCDischarge} and
\ref{figSCDischargeLog} are similar, 
but numerical stability of the results can be different
in general case.

We also want to note that instead of the basis $Q_j(x)$
a basis $Q_j(f(x))$ can be used. Similar technique was
used in\cite{2016arXiv160204423G} for market dynamics study.
The formulas are very similar to (\ref{Cqm}) and (\ref{CDqm}),
e.g. we now have
$\Braket{df/dx Q_k(f(x))}=\sum_{l=1}^{M} Q_k(f(x_l)) (f_l-f_{l-1})\omega(x_l)$
instead of (\ref{CDqm}). The $Q_j(x)$ and $Q_j(f(x))$ bases
typically give very similar results for $n>10$,
but numerical stability is often better for  $Q_j(x)$ basis.

\subsection{\label{Tracker} Supercapacitor With Track-Etched Membrane}
\begin{figure}
  \includegraphics[width=4cm]{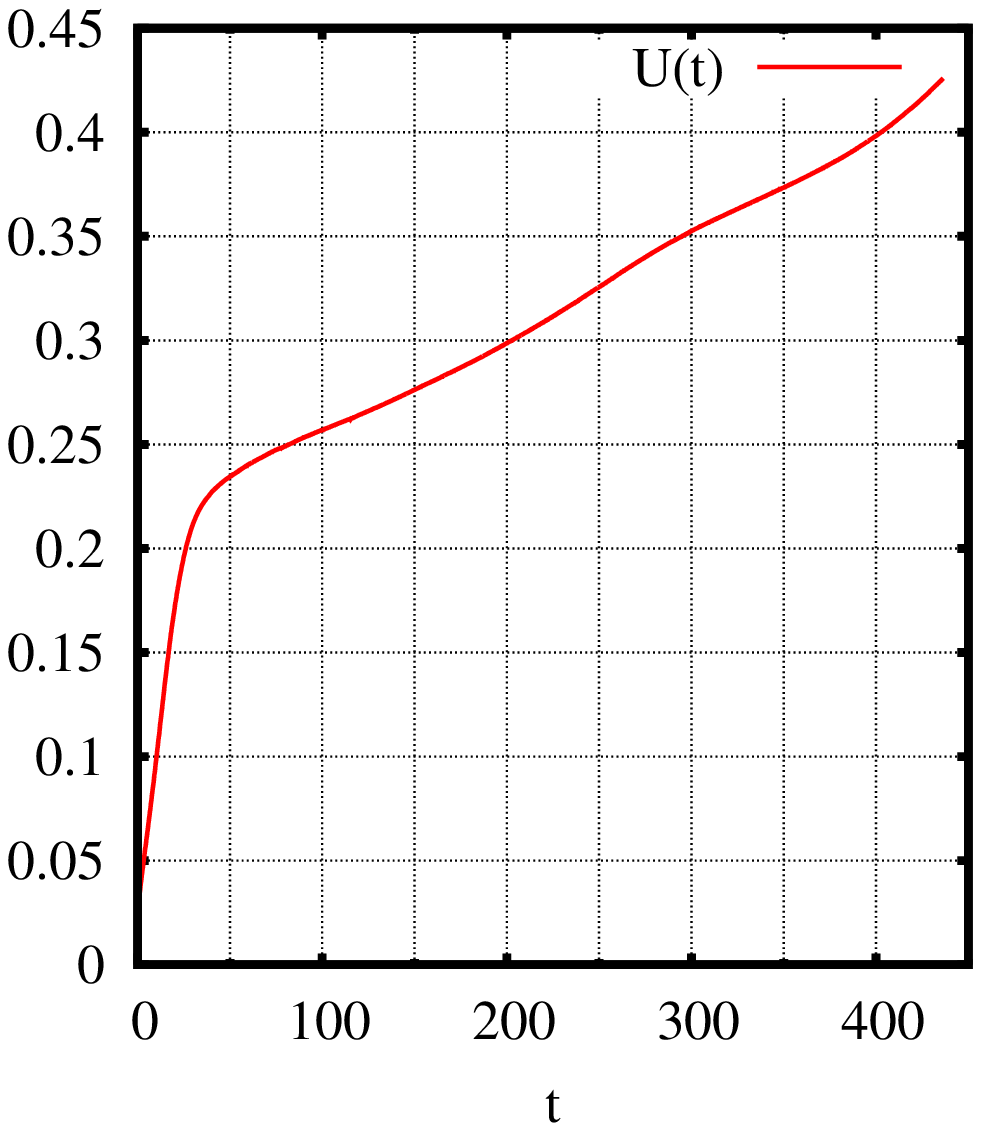}
  \includegraphics[width=4cm]{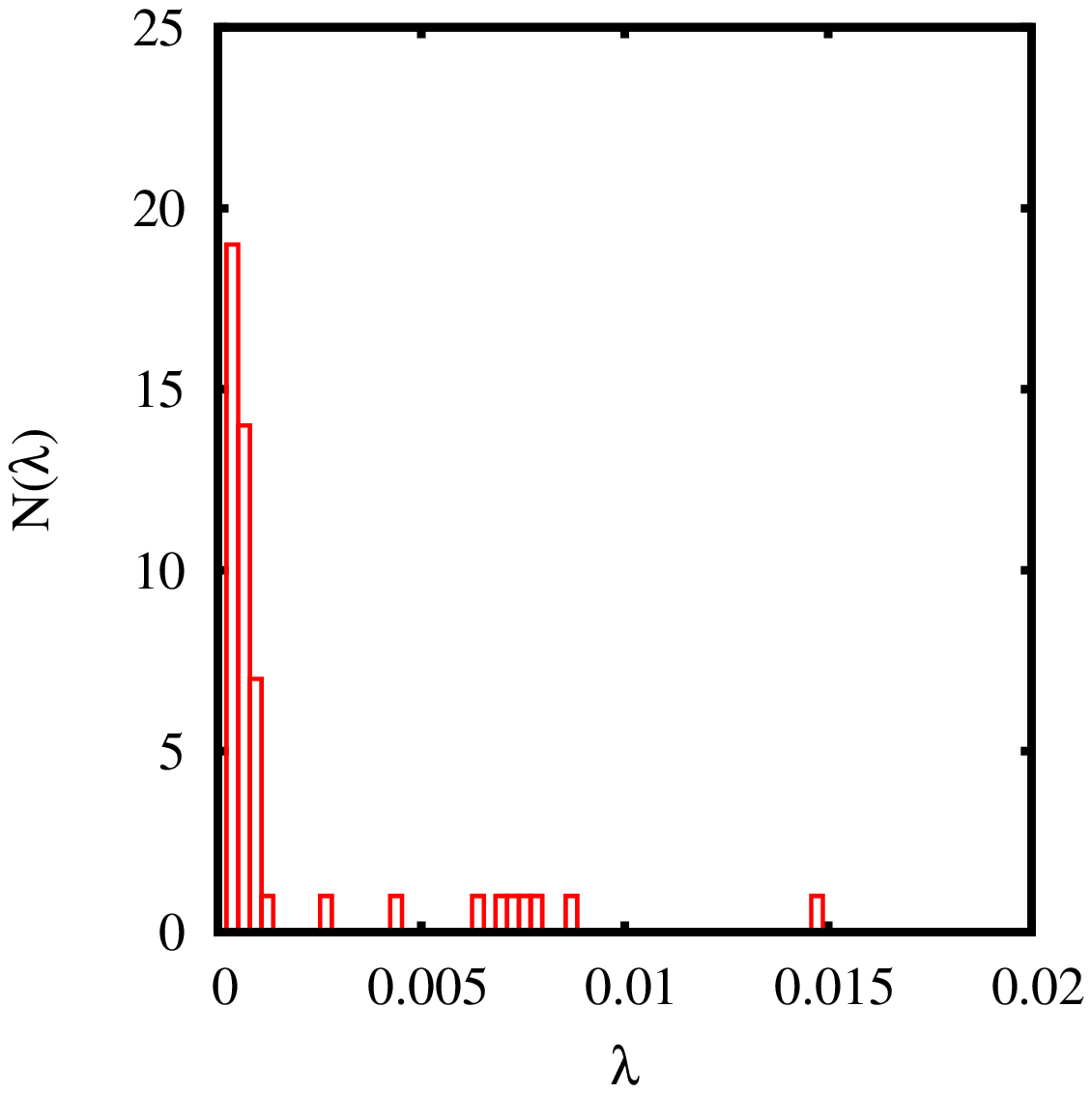}
  \includegraphics[width=9cm]{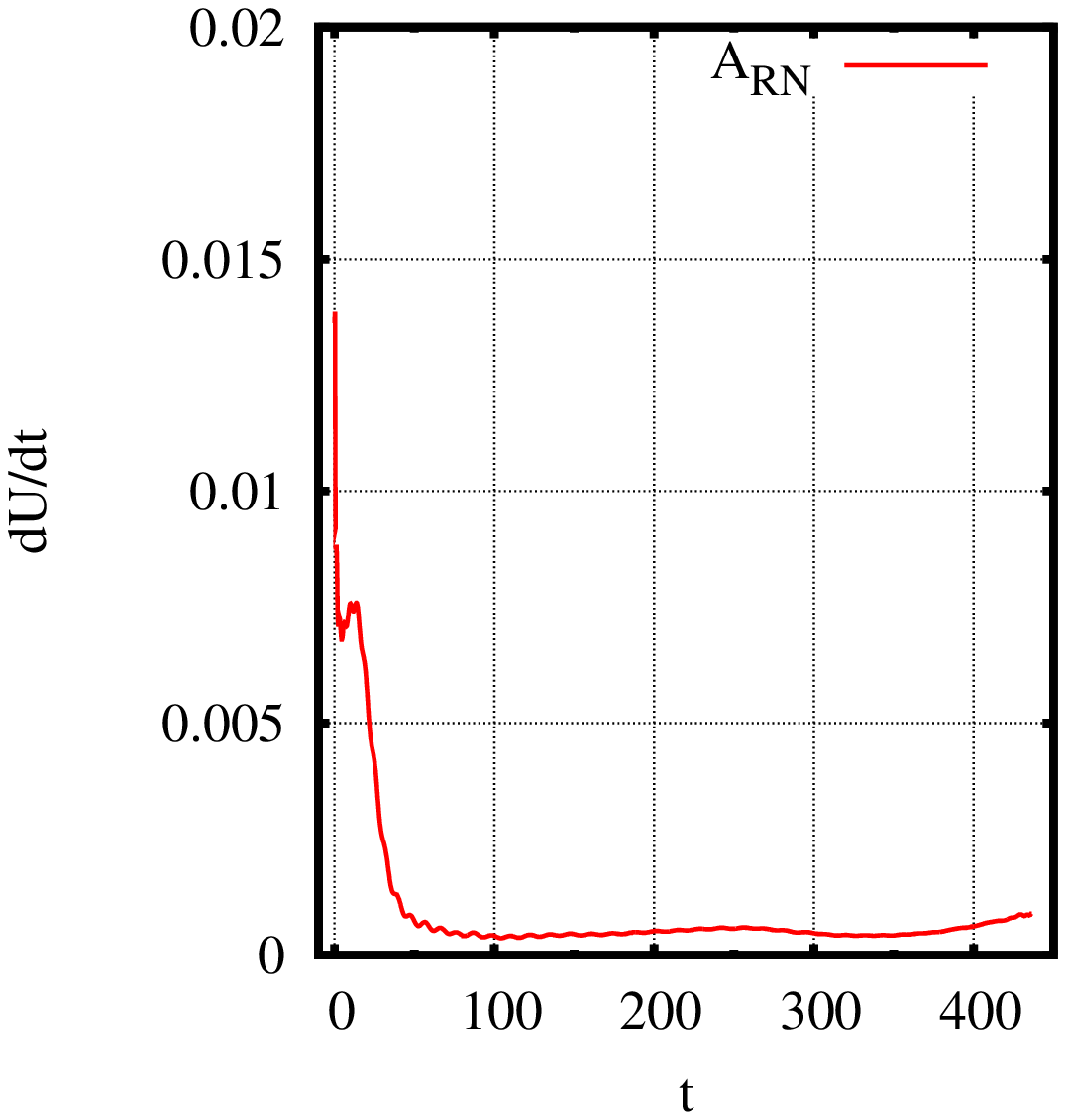}
  
  \caption{\label{figSCTrackEtchedCharge}
    Charging process of 
  a supercapacitor with track-etched membrane.
  Potential $U$ as a function of time (top left).
  $dU/dt$ (bottom) calculated as (\ref{RNsimple}) with $n=50$,
  $g=dU/dt$, $x=t$ and
  (\ref{CDqm}) moments.  
  The distribution of $\lambda$ is calculated with $n=50$ (top right).
  }
\end{figure}
In Fig. \ref{figSCTrackEtchedCharge}
 charging process data $U(t)$ with $I=2.5\cdot 10^{-5}$A
for supercapacitor
with  superionic solid state electrolyte $RbAg_4I_5$
and track-etched membrane is presented.
(Because electrolyte stability window of $RbAg_4I_5$ is
$0.55-0.60$V, maximal charging
potential is limited to $0.43$V. Charging time $436$sec.)
This system can undergo multiple charge--discharge
cycles and the measurement process can be completed under ten minutes.
For these reasons this system is extemely 
convenient for relaxation type of signals study.
Typical for supercapacitors two--stage $U(t)$ is seen on the chart.
$dU/dt$ interpolation (using localized $\psi(x)$ from (\ref{psiyX}))
is calculated according to (\ref{RNsimple}) with (\ref{CDqm}) moments.
The most interesting is  $dU/dt$ distribution chart,
where the two intervals of different $\lambda$ are clearly seen.

\subsection{\label{StockP} Stock Price Changes}
\begin{figure}
  \includegraphics[width=7.5cm]{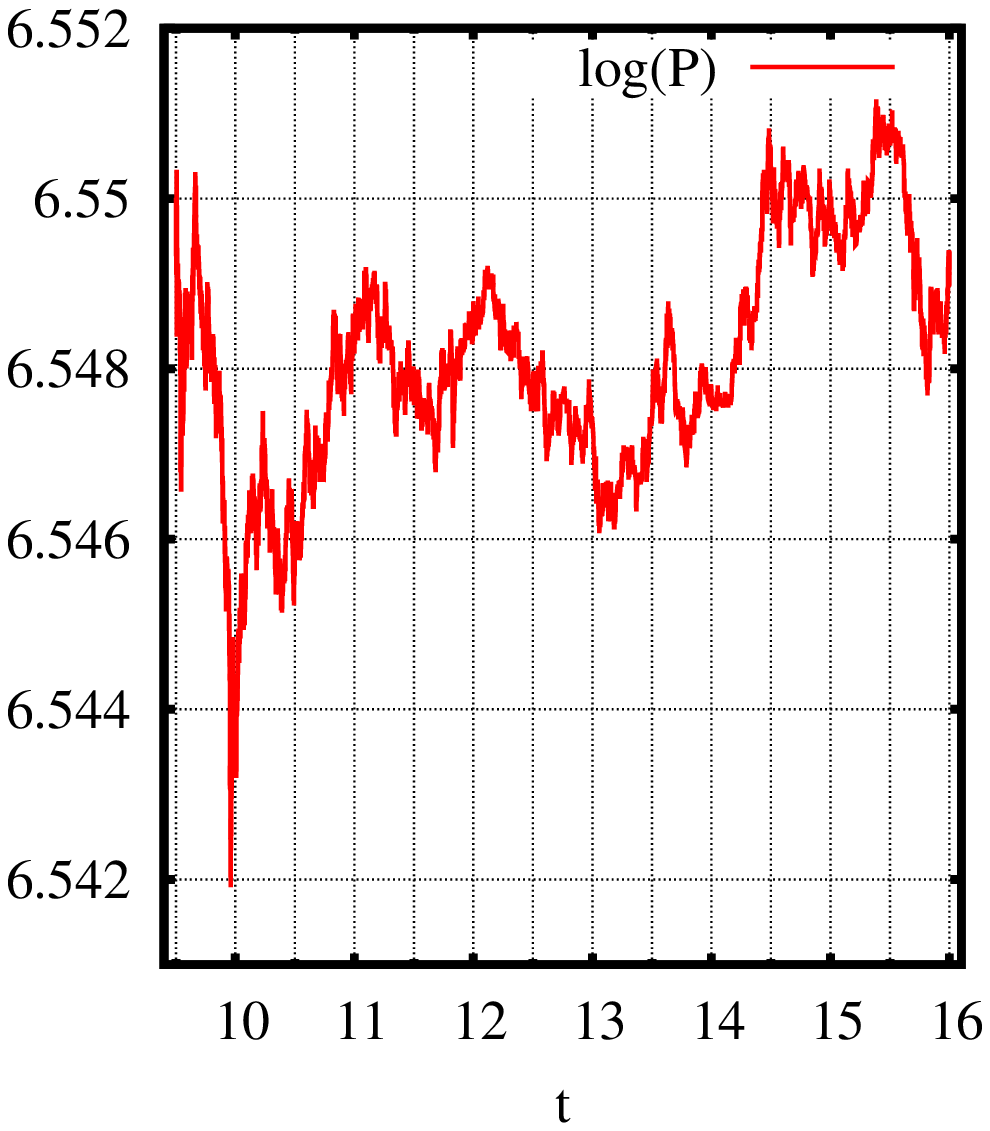}
  \includegraphics[width=7.5cm]{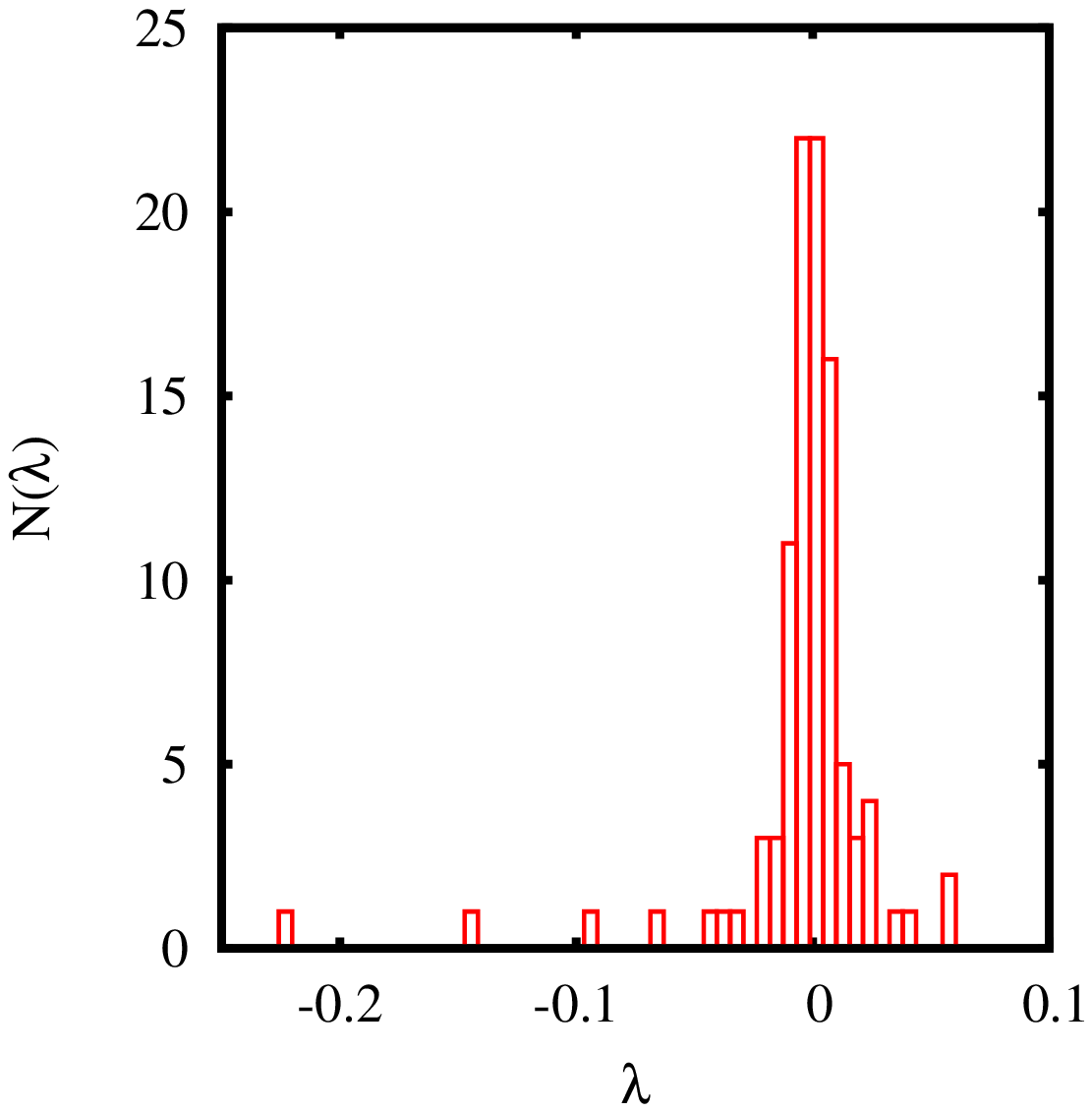}
  \caption{\label{stockF}
    The logarithm of AAPL stock price on September, 20, 2012.
    $\lambda$ distribution of $d\ln P /dt$,
    calculated with $n=100$.
  }
 \end{figure} 

As we noted above the technique can be applies to arbitrary $f$ and $x$.
Consider AAPL stock price on
September, 20, 2012, same day we used in \cite{2015arXiv151005510G,2016arXiv160204423G}
and obtain the distribution of $d(\ln P))/dx$, (here $x$ is time in hours).
The result is presented in Fig. \ref{stockF} for $n=100$.
A very important feature of this approach is that
it can be successfully applied to non--Gaussian distributions.
The distribution of returns in Fig. \ref{stockF} is clearly non--Gaussian,
but, in contrast with all the $L^2$ approaches,
matrix spectrum density analysis
allows to estimate probability distribution
even in a strongly non--Gaussian case\cite{malyshkin2015norm}.
To reproduce the result presented in Fig. \ref{stockF} numerically see appendix \ref{PriceCalc} below.

\section{\label{Applicab} Radon--Nikodym Approach Applicability}

There are several important issues, affecting results quality.
\begin{itemize}
\item Numerical instability. The problem arise
  typically for $n>4$ in monomials basis
  when calculating the moments (\ref{Cqm}) and  (\ref{CDqm}),
  and then calculating $\Braket{Q_j g Q_k}$ from $\Braket{g Q_k}$ moments,
  what requires stable multiplication (\ref{cmul}) of basis functions.
  With a proper basis choice (like Chebyshev or Legendre)
  the calculation are stable\cite{laurie1979computation,beckermann1996numerical,2015arXiv151005510G,2015arXiv151101887G} for $n$ up to 100-150
 for 64bit double precision computer arithmetic.
  Stable calculations at even higher $n$ can
  be achieved with a transition to high--precision
  arithmetic.
\item Data quality. While numerical instability in
  (\ref{Cqm}) and  (\ref{CDqm}) can be overcome,
  missed data in the sample strongly affect the results.
\item
  In situations, when high order derivatives are of interest,
  a special care should be taken of the boundaries.
  For example, $\Braket{d^2f/dx^2 Q_j}$ moments can be easily obtained
  from $\Braket{df/dx Q_j}$ moments (\ref{CDqm}) using integration
  by parts for all practically interesting $\omega(x)$.
  However, boundary terms strongly affect $\Braket{d^2f/dx^2 Q_j}$ moments
  what make  $d^2f/dx^2$ distribution much more difficult to study.
\item
  Distribution results are of non--interpolatory type.
  Radon--Nikodym approach to $df/dx$ distribution estimation
  finds the state of virtual Hamiltonian, such that: 
  for any $\lambda^{[i]}$ in $N(\lambda)$ distribution
  there is exist a probability density $(\psi^{[i]}(x))^2$, such
  that $\lambda^{[i]}=\Braket{df/dx (\psi^{[i]}(x))^2}/\Braket{(\psi^{[i]}(x))^2}$.
  This is different from often used results of interpolatory type,
  where the answer cannot be interpreted as averaging
  with always positive weight.
  In this sense even a little peak in $N(\lambda)$ is important
  as it correspond to \textsl{actual averaging with positive weights}.
\end{itemize}

\section{\label{discuss} Discussion}
Generalized Radon--Nikodym approach,
where the probability density is built first
and then the value of dynamic characteristic is obtained by
averaging with this probability density
is developed and implemented numerically.
The  approach
is most advantageous to relaxation dynamics study,
because the spectrum of relaxation time can be directly obtained
as the spectrum of specially constructed operator (virtual Hamiltonian).
This is especially important, because 
Fourier or Wavelet analysis are poorly applicable
to  relaxation  dynamics. Laplace analysis
is not applicable because of 
sample insufficient size and discretization noise.
In contrast with Laplace approach, Radon--Nikodym
can be effectively applied (and processed numerically) to sampled data.
Despite ideological differences with Fourier
(operator spectrum vs. $L^2$ norm)
generalized Radon--Nikodym approach is applicable to relaxation processes
as widely as Fourier approach is applicable to oscillatory processes.
Software product, implementing the theory is provided.

\appendices
\section{\label{appnum} Computer Code Numerical Calculations}
This software product is licensed under GPL version 3 license
and can be downloaded from   \cite{polynomialcode} website.
If for integration to a commercial product different
license is required the authors are open to
consideration of such a request.

The program read timeserie pairs $(x_l,f_l)$, calculate
$\Braket{Q_j Q_k}$, $\Braket{Q_j f Q_k}$ and $\Braket{Q_j df/dx Q_k}$ matrices
according to (\ref{qm}), (\ref{Cqm}) and (\ref{CDqm})
with $\omega(x)=1$.
Then the following calculations are performed:
\begin{itemize}
\item Least squares interpolation (\ref{leastsq}).
\item Radon--Nikodym interpolation (\ref{RNsimple}).
\item The spectrum of (\ref{GEV})  with $M^{L}_{jk}=\Braket{ Q_j f Q_k}$ and $M^{R}_{jk}=\Braket{ Q_j Q_k}$.
\item The spectrum of (\ref{GEV})  with $M^{L}_{jk}=\Braket{ Q_j df/dx Q_k}$ and $M^{R}_{jk}=\Braket{ Q_j Q_k}$.
\item The spectrum of (\ref{GEV})  with $M^{L}_{jk}=\Braket{ Q_j df/dx Q_k}$ and $M^{R}_{jk}=\Braket{ Q_j f Q_k}$.
\end{itemize}

The $(x_l,f_l)$ timeserie data must input as
two column tab--separated file, the lines starting with ``$|$''
are considered to be comments. Assume the $(x_l,f_l)$ input is saved to
\verb+input_file.dat+. The program
\begin{verbatim}
java com/polytechnik/algorithms/\
ExampleRadonNikodym_F_and_DF \
   input_file.dat n flagDX
\end{verbatim}
needs three arguments to be specified.
\begin{itemize}
\item \verb+input_file.dat+: input $(x_l,f_l)$ timeserie file
\item \verb+n+: basis dimension, typical value is between 4 and 100.
\item \verb+flagDX+: either \verb+sampleDX+ or \verb+analyticalDX+.
  Depending on its value the (\ref{qm}) is calculated
  either numerically or analytically.
  The result typically does not depend on this parameter
 unless input timeserie has too few observations.
\end{itemize}

The program outputs six files (names are hardcoded):
\begin{itemize}
\item \verb+RN_interpolated.dat+  
  Interpolation result files. Tab--separated file,
  the columns correspond to the following:
  \begin{enumerate}
  \item \verb+x+ original $x_l$.
  \item \verb+f_orig+ original $f_l$.
  \item \verb+f_RN+ interpolation $f(x_l)$ according to (\ref{RNsimple}).
  \item \verb+f_LS+ interpolation $f(x_l)$(\ref{leastsq}).
  \item \verb+df_RN+  interpolation $df/dx(x_l)$ according to (\ref{RNsimple}).
  \item \verb+df_LS+ interpolation $df/dx(x_l)$ according to (\ref{leastsq}).
  \item \verb+df_RN_byparts+ and \verb+df_LS_byparts+ are similar to
    \verb+df_RN+ and \verb+df_LS+
    but the $\Braket{df/dxQ_k}$ moments
    are calculated from the $\Braket{fQ_k}$  moments using integration by parts,
    this can be easily done in case of constant weight $\omega(x)=1$.
  \end{enumerate}
\item  \verb+EV_RN_interpolated.dat+ has the same results
  as \verb+RN_interpolated.dat+, but the formulas (\ref{RNsimpleEVBasis}) and (\ref{leastsqEVBasis})
  are used instead of (\ref{RNsimple}) and (\ref{leastsq}).
  The results should be identical within computer real numbers arithmetic precision.
  Different results indicate numerical instability.
\item \verb+QQf_QQ_spectrum.dat+, \verb+QQdf_QQ_spectrum.dat+,
  \verb+QQdfbyparts_QQ_spectrum.dat+
  and \verb+QQdf_QQf_spectrum.dat+. These are the files with eigenvalues
  of generalized eigenvalues (\ref{GEV}) problems
  solved for the following $\Big( M^{L}_{jk} ; M^{R}_{jk} \Big)$
  pairs: $\Big(\Braket{ Q_j f Q_k} ; \Braket{ Q_j Q_k}\Big)$,
  $\Big(\Braket{ Q_j df/dx Q_k} ; \Braket{ Q_j Q_k}\Big)$ (calculate directly),
  $\Big(\Braket{ Q_j df/dx Q_k} ; \Braket{ Q_j Q_k}\Big)$ (calculate from $\Braket{ Q_j f}$ using integration by parts)
  and $\Big(\Braket{ Q_j df/dx Q_k} ; \Braket{ Q_j f Q_k}\Big)$
  respectively.
  The files have three tab--separated columns: index $(0..n-1)$
  eigenvalue $\lambda^{[i]}$ and the value of
  $x^{[i]}_{est}=\Braket{\left(\psi^{[i]}(x)\right)^2 x}\Big/\Braket{\left(\psi^{[i]}(x)\right)^2}$
  estimated as in (\ref{xest}).
  The eigenvalues are sorted in ascending order.
  If right--side (\ref{GEV}) matrix $M^{R}_{jk}$ is not positively defined
  all eigenvalues are set to NaN.
\end{itemize}
All the calculations are performed in Chebyshev
polynomial basis as providing a very good numerical stability
and having numerically always  stable basis multiplication (\ref{cmul}).
The bases of Legendre, Laguerre, Hermite and monomials
are also implemented. The results are invariant with respect
to basis selection, but numerical stability can be drastically different.
For this reason the file com/ polytechnik/ algorithms/ ExampleRadonNikodym\_F\_and\_DF.java
may be modified
\begin{verbatim}
final OrthogonalPolynomialsABasis 
//Q=new OrthogonalPolynomialsChebyshevBasis()
Q=new OrthogonalPolynomialsLegendreBasis()
//Q=new OrthogonalPolynomialsMonomialsBasis()
;
\end{verbatim}
to check the stability of the results in Legendre basis,
that also provide a good stability in most cases.
The results should be identical within computer real numbers arithmetic precision.

\subsection{\label{install} Software Installation And Testing}

\begin{itemize}
\item Install java 1.8 or later.
\item Download the file Electrochemistry.zip from \cite{polynomialcode} website.
\item Decompress and recompile the program
\begin{verbatim}
unzip Electrochemistry.zip
javac -g com/polytechnik/*/*java
\end{verbatim}
\item Run the test, corresponding to the results of \cite{LGChemLGHG2}
  data. Obtain the data presented in Fig. \ref{figLiIonLG}.
\begin{verbatim}
java com/polytechnik/algorithms/\
ExampleRadonNikodym_F_and_DF \
   echem/lg10.dat 50 sampleDX
\end{verbatim}
The six files \\
\verb+RN_interpolated.dat+, 
\verb+EV_RN_interpolated.dat+, 
\verb+QQf_QQ_spectrum.dat+, 
\verb+QQdf_QQ_spectrum.dat+, 
\verb+QQdfbyparts_QQ_spectrum.dat+
\verb+QQdf_QQf_spectrum.dat+
have to match identically to the files provided in Electrochemistry.zip: \\
\verb+echem/lg10_RN_interpolated.dat+, 
\verb+echem/lg10_EV_RN_interpolated.dat+, 
\verb+echem/lg10_QQf_QQ_spectrum.dat+, 
\verb+echem/lg10_QQdf_QQ_spectrum.dat+,
\verb+echem/lg10_QQdfbyparts_QQ_spectrum.dat+
\verb+echem/lg10_QQdf_QQf_spectrum.dat+.

\end{itemize}

\subsection{\label{PriceCalc} Stock Price Change Rate Distribution Calculation}
To obtain the data presented in Fig. \ref{stockF} follow these steps
\begin{itemize}
\item Download from \cite{polynomialcode}  data file S092012-v41.txt.gz ,
  the  same one used in Ref. \cite{2016arXiv160204423G}.
\item Extract from the file S092012-v41.txt.gz AAPL executed trades
  and save the data to aapp.csv same as \cite{2016arXiv160204423G}.
\begin{verbatim}
java com/polytechnik/itch/DumpDataTrader \
   S092012-v41.txt.gz AAPL >aapl.csv
\end{verbatim}
\item Convert the data from aapl.csv to market hours (9:30 to 16:00)
  and the form  (time, $\ln P$), then run Radon--Nikodym code example from previous section.
\begin{verbatim}
java  echem/Extract23cols \
   aapl.csv >aapl2cols.csv
java com/polytechnik/algorithms/\
ExampleRadonNikodym_F_and_DF \
   aapl2cols.csv 100 sampleDX
\end{verbatim}

\end{itemize}

\section*{Acknowledgment}
  Vladislav Malyshkin would like to thank
  N. S. Averkiev and V. Yu. Kachorovskii
  for fruitful discussions of
  quantum mechanical description of classical experiments.

\bibliographystyle{IEEEtran}
\bibliography{LD}

\end{document}